\documentclass[11pt]{article}

\usepackage{amsmath,amssymb,amsthm,amsfonts}
\usepackage{bm}
\usepackage{geometry}
\usepackage{graphicx}
\usepackage{multirow}
\usepackage{caption}
\usepackage{subcaption}
\usepackage{authblk}   
\usepackage{placeins}
\usepackage{tikz}
\usepackage{microtype}

\newtheorem{theorem}{Theorem}
\newtheorem{lemma}{Lemma}
\newtheorem{remark}{Remark}
\newtheorem{assumption}{Assumption}

\title{Multicontinuum Poroelasticity Model. Analysis and Computation.}

\author[1]{Jing Zhang}
\author[2]{Mohammed Al Kobaisi}
\author[2]{Dmitry Ammosov}
\author[1]{Jian Huang}
\author[1]{Yunqing Huang}

\affil[1]{School of Mathematics and Computational Science, Xiangtan University \& National Center for Applied Mathematics in Hunan \& Hunan Key Laboratory for Computation and Simulation in Science and Engineering \& Hunan Research Center of the Basic Discipline Fundamental Algorithmic Theory and Novel Computational Methods, Xiangtan, Hunan Province, PRC. Email: jingzhang@xtu.edu.cn (Jing Zhang).}
\affil[2]{Chemical and Petroleum Engineering Department, Khalifa University of Science and Technology, Abu Dhabi, 127788, UAE.}

\date{}

\begin{document}

\maketitle

\begin{abstract}
In this work, we develop and analyze a macroscopic multicontinuum model for Biot poroelasticity. The full macroscopic multicontinuum system is derived by projecting the microscopic weak formulation onto continuum-dependent Galerkin spaces enriched with first-order multiscale corrections. This construction retains generalized pressure storage and coupling terms, with effective coefficients inherited from the microscopic bilinear forms. Under suitable coefficient assumptions, we establish well-posedness and energy stability of the macroscopic multicontinuum model. A fully discrete scheme combining the Galerkin finite element method in space with the backward Euler method in time is then developed. Its well-posedness and stability are proved, and a priori error estimates are derived using appropriate projection operators and energy arguments. Finally, numerical experiments for periodic and non-periodic high-contrast media demonstrate the accuracy of the multicontinuum model with respect to the averaged fine-scale solution.
\end{abstract}


\section{Introduction}
\label{sec:introduction}

Poroelasticity describes the interaction between fluid flow and
solid deformation in porous media and plays an important role in
geomechanics, subsurface engineering, and biological tissue modeling.
The classical Biot model couples mechanical equilibrium with fluid
mass conservation through the dependence of stress on pore pressure
and the contribution of volumetric deformation to fluid storage
\cite{Biot1941,Biot1955,Biot1956}. In heterogeneous porous media \cite{frias2004stochastic}, the elastic,
hydraulic, and coupling coefficients may vary over multiple spatial
scales. Resolving these variations directly can be computationally
expensive, particularly when strong material contrasts and complex
microstructures are present. Meanwhile,
the coupling between fluid flow and solid deformation complicates
both the construction of effective models and the analysis of
their well-posedness and stability. This motivates the development of
macroscopic models that retain the essential effects of fine-scale
heterogeneity while reducing computational complexity.

To address the multiscale nature of heterogeneous media, various multiscale numerical methods have been developed. Multiscale finite element methods construct local basis functions by solving fine-scale problems on coarse-grid blocks, thereby incorporating unresolved heterogeneities into a reduced coarse-grid approximation. The generalized multiscale finite element method (GMsFEM) further employs snapshot spaces and local spectral problems to systematically select the dominant modes, and has been applied to linear and nonlinear poroelasticity problems \cite{EfendievHou2009,Brown2016,Brown2016Nonlinear}. 
The constraint energy minimizing generalized multiscale finite
element method (CEM-GMsFEM) further constructs basis functions
through constrained energy minimization on oversampled regions
\cite{ChungEtAl2018CEM}.
It has also been applied to poroelasticity problems with high-contrast coefficients \cite{FuEtAl2019}. 
Although these approaches can substantially
reduce the computational cost of resolving fine-scale
heterogeneity, their direct output is typically a discrete
coarse-grid system rather than a continuous macroscopic model expressed in terms of explicit governing laws.

Multicontinuum modeling offers an alternative description for heterogeneous porous media in which a single macroscopic pressure variable may be insufficient to represent distinct fluid responses. A representative example is the multiple-network poroelastic theory (MPET), which extends the classical Biot model to a deformable porous skeleton permeated by several interacting fluid networks \cite{wilson1982theory,LeeEtAl2019}. Building on double-porosity
and multiple-porosity descriptions, MPET provides a flexible
framework for representing coupled flow and deformation in
complex porous systems and has been applied to biological
tissues, particularly the human brain \cite{GuoEtAl2019,guo2020multiple}.
However, the predictive capability of these models depends
on the specification of effective storage, transfer, and
poromechanical coupling coefficients. Determining these
macroscopic coefficients, especially cross-porosity storage
coefficients, can be challenging because experimental data
capable of quantifying them are often limited
\cite{guo2020multiple,Ammosov2025}.
This motivates
systematic approaches that relate macroscopic multicontinuum
parameters to the underlying microscopic properties.

The multicontinuum homogenization method provides a systematic framework for deriving such macroscopic multicontinuum models from local microscopic problems. Its central idea is to introduce several macroscopic variables, or continua, at each coarse-scale location and represent the fine-scale solution through continuum-dependent basis functions obtained from constrained local problems. The constraints are associated with prescribed continuum averages, while the local problems are posed on representative volume elements and, in particular, on oversampled regions to reduce the influence of artificial boundary conditions \cite{EfendievLeung2023}. This framework has been further developed to treat general systems, different types of constraints, and mixed or nonsymmetric formulations \cite{ChungEtAl2024,ammosov2026multicontinuum,al2026dynamic}. 
For poroelasticity,
Ammosov et al.~\cite{Ammosov2025} constructed coupled
displacement-pressure cell problems, derived generalized
multicontinuum models, and numerically investigated both
full and simplified models. Building on these developments,
a systematic analysis of the full macroscopic system is
needed to characterize the role of generalized storage and
coupling terms in its variational structure, well-posedness,
and energy stability.

Finite element methods for Biot poroelasticity have been studied extensively. Phillips and Wheeler developed mixed and continuous Galerkin formulations for the classical Biot system and established convergence and error estimates for continuous- and discrete-in-time approximations \cite{PhillipsWheeler2007a,PhillipsWheeler2007b}. For multiple-network poroelasticity, Lee et al.~\cite{LeeEtAl2019} proposed and analyzed a mixed finite element formulation, including energy estimates and a priori error estimates that are robust in nearly incompressible and low-storage regimes. These works provide important analytical foundations for the numerical treatment of coupled poromechanical systems. The full multicontinuum model considered here, however, contains additional generalized pressure-storage terms and gradient-dependent coupling terms generated by first-order multiscale corrections. Its analysis therefore requires suitable boundedness and coercivity conditions on the assembled bilinear forms, together with a careful treatment of the reciprocity relations between the coupling operators. In particular, the cancellation of reciprocal coupling terms in the continuous energy estimate must be preserved by the discretization to establish corresponding stability and error bounds.

The purpose of this paper is to investigate the well-posedness
and finite element approximation of the full macroscopic
multicontinuum Biot system.
Starting from the microscopic Biot weak
formulation, we construct continuum-dependent Galerkin trial
and test spaces enriched with first-order multiscale corrections.
Projection onto these spaces defines the effective macroscopic
bilinear forms and makes their symmetry and coupling reciprocity
explicit. Under suitable boundedness, coercivity, and
reciprocity assumptions, we establish the existence and
uniqueness of weak solutions and derive an energy stability
estimate. We then introduce a Galerkin finite element discretization in space together with the backward Euler method in time. For the resulting fully discrete scheme,
we prove well-posedness and energy stability and derive
a priori error estimates using appropriate projection operators
and energy arguments. 
Numerical experiments for periodic and non-periodic high-contrast media demonstrate the accuracy of the multicontinuum model with respect to the averaged fine-scale solution.

The remainder of the paper is structured as follows. 
Section \ref{sec:model} introduces the microscopic Biot problem and the variational formulation of the full multicontinuum model, and describes its fully discrete finite element scheme. 
Section \ref{sec:wellposedness} proves the existence, energy stability of the continuous model, and the well-posedness and stability of the discrete scheme. 
Section \ref{sec:error} derives a priori error estimates. Section \ref{sec:numerics} reports numerical experiments. Finally, Section~\ref{sec:conclusions} summarizes the main
findings.

Through out the paper we use $C$, with or without subscript, to denote a positive constant, which could have different values at different
appearances.

\section{Preliminaries. Multicontinuum models}
\label{sec:model}
We derive the multicontinuum macroscopic equations directly from the
variational formulation of the Biot model using the multicontinuum homogenization method.
\subsection{Biot model and variational formulation}
Let $\Omega\subset\mathbb{R}^d$ be a bounded domain and let
$J=(0,T]$ denote the time interval, where $T>0$ is the final time. We consider the quasi-static
Biot system
\begin{equation}
\begin{aligned}
&-\nabla\cdot \sigma( u )
+\alpha\nabla p
= f ,
&&\quad \text{in }\Omega\times J,\\
&c_0 \partial_t p+\alpha \nabla\cdot\partial_t u 
-\nabla\cdot\left(k\nabla p\right)
=g,
&&\quad \text{in }\Omega\times J.
\end{aligned}
\label{eq:biot_model}
\end{equation}
Here, $ u$ denotes the displacement, $p$ is the pore
pressure,
and $\sigma( u ) = 2\mu \varepsilon( u )+\lambda(\nabla\cdot u) I$ denotes the stress
tensor,
$
\varepsilon( u )
=
\frac12\left(\nabla u +\nabla u ^{T}\right)
$
denotes the strain tensor, $\alpha$ is the Biot-Willis coefficient, $c_0$ is
the constrained storage coefficient, $k$ is the permeability tensor,
and $ f $ and $g$ denote the prescribed body force and fluid source,
respectively.
For simplicity, we impose homogeneous Dirichlet boundary conditions,
\begin{equation*}
u= 0,
\qquad
p=0,
\qquad
\text{on }\partial\Omega\times J,
\label{eq:biot_bc}
\end{equation*}
together with the initial condition
\begin{equation*}
p( x ,0)=p_0( x ),
\qquad
 x \in\Omega.
\label{eq:biot_ic}
\end{equation*}
Define
\[
W_u=[H_0^1(\Omega)]^d,
\qquad
W_p=H_0^1(\Omega).
\]
The weak formulation reads: find $( u,p)\in W_u\times W_p$ such that,
for every $t\in J$,
\begin{equation}
\begin{aligned}
&a_u( u, v)-b(p, v)=l_u( v),&&  v\in W_u, \\
&c(\partial_t p,q)+b(q,\partial_t u)+a_p(p,q)=l_p(q),&& q\in W_p,
\end{aligned}
\label{eq:weakbiot}
\end{equation}
where
\begin{align*}
&a_u( u, v)
=
\int_\Omega
\sigma( u):
\varepsilon(\bm{v}),&& b(p, v)
=
\int_\Omega
\alpha p\,\nabla\!\cdot  v,\\
&a_p(p,q)
=
\int_\Omega
k
\nabla p\cdot\nabla q, &&c(p,q)
=
\int_\Omega
c_0pq,\\
& l_u( v) = \int_\Omega  f  v, && l_p(q) = \int_\Omega gq,
\end{align*}
\subsection{Multicontinuum approximation}
This multicontinuum formulation follows \cite{Ammosov2025}; however, we split the 
expansion for displacement and pressure fields.
Let $N_c$ denote the number of continua. We introduce the
multicontinuum approximations
\begin{align}
 u_H
&=
\sum_{i=1}^{N_c}
\left(
N_i^uU_i+M_i^u\cdot\nabla U_i
\right),
\label{eq:utrial}\\
p_H
&=
\sum_{i=1}^{N_c}
\left(
N_i^pP_i+M_i^p\cdot\nabla P_i
\right).
\label{eq:ptrial}
\end{align}
where $U_i$ and $P_i$ are the macroscopic displacement and pressure
variables associated with the $i$th continuum. The functions
$N_i^u$, $M_i^u$, $N_i^p$, and $M_i^p$ denote the corresponding
multicontinuum basis functions. 
These indices are suppressed for readability, and all
products below denote the corresponding tensor
contractions.
In \cite{Ammosov2025}, the multicontinuum
expansion for the displacement has pressure
terms and vice versa. The resulting macroscopic equations are the same. Since the goal of this paper is to analyze macroscopic equations, for
simplicity, we consider split multicontinuum expansion (\ref{eq:utrial}) and (\ref{eq:ptrial}).
The test functions are selected from the same approximation
spaces:
\begin{align}
 v_H
&=
\sum_{j=1}^{N_c}
\left(
N_j^uV_j+M_j^u\cdot\nabla V_j
\right),
\label{eq:vtest}\\
q_H
&=
\sum_{j=1}^{N_c}
\left(
N_j^pQ_j+M_j^p\cdot\nabla Q_j
\right).
\label{eq:qtest}
\end{align}
The construction of the effective coefficients is performed locally on
each Representative Volume Element (RVE) region $R_{\omega}$. Within $R_{\omega}$, the macroscopic variables
and their first-order spatial derivatives are regarded as prescribed
parameters, while the microscopic differential operators act on the
corresponding local basis functions. Thus, the resulting expressions
should be interpreted as a local first-order approximation of the
multicontinuum representation, rather than as an exact substitution of
spatially varying macroscopic fields into the microscopic differential
operators.
For a microscopic bilinear form $a(\cdot,\cdot)$, we denote its
normalized restriction to $R_{\omega}$ by
\begin{equation*}
a^{R_{\omega}}(\phi,\psi)
:=
\frac{1}{|R_{\omega}|}
a(\phi,\psi;R_{\omega}),
\label{eq:local_bilinear_form}
\end{equation*}
where $a(\phi,\psi;R_{\omega})$ denotes the contribution of the region $R_{\omega}$ to
the microscopic bilinear form. The other microscopic bilinear forms are defined similarly.
The effective coefficients are then
computed locally on each RVE region and subsequently
interpreted as spatially varying macroscopic coefficient fields.

\subsection{Effective operators}
\label{subsec:effective}
Substitution of the multicontinuum trial and test functions \eqref{eq:utrial}-\eqref{eq:qtest} into
\eqref{eq:weakbiot} yields a macroscopic system for $(U,P)$. We first
give the effective coefficient blocks generated by the microscopic
bilinear forms.
The effective elastic bilinear form is
\begin{equation*}
\begin{aligned}
a_u( u_H, v_H)
=
\sum_{i,j=1}^{N_c}\int_\Omega
\Big[
&
A^{uu}_{ij}\nabla U_i\cdot\nabla V_j
+
B^{uu}_{ij}\nabla U_i\cdot V_j\\
&
+
\bar B^{uu}_{ij}U_i\cdot\nabla V_j
+
C^{uu}_{ij}U_i\cdot V_j
\Big]\,d x,
\end{aligned}
\label{eq:elastic_expansion}
\end{equation*}
where
\begin{equation*}
\begin{aligned}
A^{uu}_{ij}&=a_u^{R_{\omega}}(M_i^u,M_j^u),
&
B^{uu}_{ij}&=a_u^{R_{\omega}}(M_i^u,N_j^u),\\
\bar B^{uu}_{ij}&=a_u^{R_{\omega}}(N_i^u,M_j^u),
&
C^{uu}_{ij}&=a_u^{R_{\omega}}(N_i^u,N_j^u).
\end{aligned}
\label{eq:elastic_coefficients}
\end{equation*}
Evaluations involving tensor-valued basis functions
are understood componentwise.
Similarly, the effective diffusion form is
\begin{equation*}
\begin{aligned}
a_p(p_H,q_H)
=
\sum_{i,j=1}^{N_c}\int_\Omega
\Big[
&
A^{pp}_{ij}\nabla P_i\cdot\nabla Q_j
+
B^{pp}_{ij}\cdot\nabla P_i\,Q_j\\
&
+
\bar B^{pp}_{ij}P_i\cdot\nabla Q_j
+
C^{pp}_{ij}P_iQ_j
\Big]\,d x,
\end{aligned}
\label{eq:effective_diffusion}
\end{equation*}
with
\begin{equation*}
\begin{aligned}
A^{pp}_{ij}&=a_p^{R_{\omega}}(M_i^p,M_j^p),
&
B^{pp}_{ij}&=a_p^{R_{\omega}}(M_i^p,N_j^p),\\
\bar B^{pp}_{ij}&=a_p^{R_{\omega}}(N_i^p,M_j^p),
&
C^{pp}_{ij}&=a_p^{R_{\omega}}(N_i^p,N_j^p).
\end{aligned}
\label{eq:diffusion_coefficients}
\end{equation*}
The effective storage form is
\begin{equation*}
\begin{aligned}
c(p_{H,t},q_H)
=
\sum_{i,j=1}^{N_c}\int_\Omega
\Big[
&
H^{pp}_{ij}\partial_t P_iQ_j
+
G^{pp}_{ij}\cdot\nabla \partial_t P_i\,Q_j\\
&
+
\bar G^{pp}_{ij}\partial_t P_i\cdot\nabla Q_j
+
D^{pp}_{ij}\nabla \partial_t P_i\cdot\nabla Q_j
\Big]\,d x,
\end{aligned}
\label{eq:effective_storage}
\end{equation*}
where
\begin{equation*}
\begin{aligned}
H^{pp}_{ij}&=c^{R_{\omega}}(N_i^p,N_j^p),
&
G^{pp}_{ij}&=c^{R_{\omega}}(M_i^p,N_j^p),\\
\bar G^{pp}_{ij}&=c^{R_{\omega}}(N_i^p,M_j^p),
&
D^{pp}_{ij}&=c^{R_{\omega}}(M_i^p,M_j^p).
\end{aligned}
\label{eq:storage_coefficients}
\end{equation*}
The coupling form is
\begin{equation*}
\begin{aligned}
b(p_H,  v_H)
=
\sum_{i,j=1}^{N_c}\int_\Omega
\Big[
&
A^{pu}_{ij}\nabla P_i\cdot\nabla V_j
+
B^{pu}_{ij}\nabla P_i\cdot V_j\\
&
+
\bar B^{pu}_{ij}P_i\cdot\nabla V_j
+
C^{pu}_{ij}P_i\cdot V_j
\Big]\,d x,
\end{aligned}
\label{eq:effective_coupling}
\end{equation*}
with
\begin{equation*}
\begin{aligned}
A^{pu}_{ij}&=b^{R_{\omega}}(M_i^p,M_j^u),
&
B^{pu}_{ij}&=b^{R_{\omega}}(M_i^p,N_j^u),\\
\bar B^{pu}_{ij}&=b^{R_{\omega}}(N_i^p,M_j^u),
&
C^{pu}_{ij}&=b^{R_{\omega}}(N_i^p,N_j^u).
\end{aligned}
\label{eq:coupling_coefficients}
\end{equation*}
The gradient-gradient coefficient $A^{pu}_{ij}$ is
included because it generally arises from the local
projection.
Using the same coupling form in the pressure equation
gives
\begin{equation*}
\begin{aligned}
b(q_H,  u_{H,t})
=
\sum_{i,j=1}^{N_c}\int_\Omega
\Big[
&
A^{up}_{ij}\nabla \partial_t U_i\cdot\nabla Q_j
+
G^{up}_{ij}\nabla \partial_t U_i\,Q_j\\
&
+
\bar G^{up}_{ij}\partial_t U_i\cdot\nabla Q_j
+
E^{up}_{ij}\partial_t U_i\,Q_j
\Big]\,d x,
\end{aligned}
\label{eq:effective_reverse_coupling}
\end{equation*}
where
\begin{equation*}
\begin{aligned}
A^{up}_{ij}&=b^{R_{\omega}}(M_j^p,M_i^u),
&
G^{up}_{ij}&=b^{R_{\omega}}(N_j^p,M_i^u),\\
\bar G^{up}_{ij}&=b^{R_{\omega}}(M_j^p,N_i^u),
&
E^{up}_{ij}&=b^{R_{\omega}}(N_j^p,N_i^u),
\end{aligned}
\label{eq:reverse_coupling_coefficients}
\end{equation*}
Similarly, the source coefficients are
\begin{equation*}
\begin{aligned}
F^{u}_{j}&=l^{R_{\omega}}(N_j^u)+l^{R_{\omega}}(M_j^u),
&
F^{p}_{j}&=l^{R_{\omega}}(N_j^p)+l^{R_{\omega}}(M_j^p).
\end{aligned}
\label{eq:source_coefficients}
\end{equation*}

The effective coefficients inherit the symmetry and reciprocity properties of the underlying microscopic bilinear forms. In particular, the elastic and diffusion operators satisfy
$$
\begin{aligned}
A^{uu}&=(A^{uu})^T,
&
C^{uu}&=(C^{uu})^T,
&
\bar B^{uu}&=(B^{uu})^T,\\
A^{pp}&=(A^{pp})^T,
&
C^{pp}&=(C^{pp})^T,
&
\bar B^{pp}&=(B^{pp})^T.
\end{aligned}
$$
Similarly, the storage coefficients satisfy
$$
H^{pp}=(H^{pp})^T,
\qquad
D^{pp}=(D^{pp})^T,
\qquad
\bar G^{pp}=(G^{pp})^T.
$$
The coupling coefficients obey the reciprocity relations
$$
\begin{aligned}
A^{up}&=(A^{pu})^T,
&
G^{up}&=(\bar B^{pu})^T,\\
\bar G^{up}&=( B^{pu})^T,
&
E^{up}&=(C^{pu})^T.
\end{aligned}
$$
Here, the transpose exchanges the corresponding trial and test component indices. These identities follow directly from the symmetry of \(a_u^{R_{\omega}}\), \(a_p^{R_{\omega}}\), and \(c^{R_{\omega}}\), together with the fact that the two sets of coupling coefficients are generated from the same microscopic coupling bilinear form \(b^{R_{\omega}}\).

\subsection{Final multicontinuum equations}
Using the effective coefficients in Section \ref{subsec:effective}, we can obtain the multicontinuum equations
\begin{equation}
\begin{aligned}
&
-\nabla\cdot(A^{uu}\nabla U)
+B^{uu}\nabla U
-\nabla\cdot(\bar B^{uu}U)
+C^{uu}U
\\
&\qquad
+\nabla\cdot(A^{pu}\nabla P)
-B^{pu}\nabla P
+\nabla\cdot(\bar B^{pu}P)
-C^{pu}P
=F_u,
\\
&
H^{pp}\partial_tP
+G^{pp}\nabla\partial_tP
-\nabla\cdot(\bar G^{pp}\partial_tP)
-\nabla\cdot(D^{pp}\nabla\partial_tP)
\\
&\qquad
-\nabla\cdot(A^{up}\nabla\partial_tU)
+G^{up}\nabla\partial_tU
-\nabla\cdot(\bar G^{up}\partial_tU)
+E^{up}\partial_tU
\\
&\qquad
-\nabla\cdot(A^{pp}\nabla P)
+B^{pp}\nabla P
-\nabla\cdot(\bar B^{pp}P)
+C^{pp}P
=F_p.
\end{aligned}
\label{eq:multicontinum_full}
\end{equation}
where $U=(U_1,\ldots,U_{N_c})^T,
\quad
P=(P_1,\ldots,P_{N_c})^T$.

We impose
\[
U=0,\qquad P=0
\quad\text{on }\partial\Omega\times J,
\qquad
P( x,0)=P_0( x),
\]
and some initial conditions.

Based on the numerical results in \cite{Ammosov2025}, we can neglect some terms,
\begin{equation*}
A^{pu}=A^{up}=0,\quad G^{pp}=\bar G^{pp}=0,\quad D^{pp} = 0.
\label{eq:coupling_truncation}
\end{equation*}
These are explicit modeling approximations and the
corresponding coefficients need not vanish in the
full local projection.
The multicontinuum equations \eqref{eq:multicontinum_full} can be written as
\begin{equation}
\begin{aligned}
&
-\nabla\cdot(A^{uu}\nabla U)
+
B^{uu}\nabla U
-
\nabla\cdot(\bar B^{uu}U)
+
C^{uu}U
\\
&\qquad
-
B^{pu}\nabla P
+
\nabla\cdot(\bar B^{pu}P)
-
C^{pu}P
=
F_u,
\\
&
H^{pp}\partial_t P
+
G^{up}\nabla \partial_t U
-
\nabla\cdot(\bar G^{up}\partial_t U)
+
E^{up}\partial_t U
\\
&\qquad
-
\nabla\cdot(A^{pp}\nabla P)
+
B^{pp}\nabla P
-
\nabla\cdot(\bar B^{pp}P)
+
C^{pp}P
=
F_p,
\end{aligned}
\label{eq:multicontinum}
\end{equation}

\begin{assumption}\label{assump:coeff_reg}
We present the Assumption as follows.
\begin{itemize}
\item[\textbf{(A1)}]
All tensor-valued coefficients are essentially bounded:
\[
\begin{aligned}
&A^{uu},B^{uu},\bar B^{uu},C^{uu},
A^{pu},B^{pu},\bar B^{pu},C^{pu},
A^{up},G^{up},\bar G^{up},E^{up},\\&
A^{pp},B^{pp},\bar B^{pp},C^{pp},
G^{pp},\bar G^{pp},D^{pp},H^{pp}
\in L^\infty(\Omega).
\end{aligned}
\]
Moreover, the effective coefficients satisfy, almost everywhere in
$\Omega$,
\[
\begin{aligned}
&A^{uu}=(A^{uu})^T,
\qquad
C^{uu}=(C^{uu})^T,
\qquad
\bar B^{uu}=(B^{uu})^T,\\
&A^{pp}=(A^{pp})^T,
\qquad
C^{pp}=(C^{pp})^T,
\qquad
\bar B^{pp}=(B^{pp})^T,\\
&H^{pp}=(H^{pp})^T,
\qquad
D^{pp}=(D^{pp})^T,
\qquad
\bar G^{pp}=(G^{pp})^T,\\
&A^{up}=(A^{pu})^T,
\qquad
G^{up}=(\bar B^{pu})^T,
\qquad
\bar G^{up}=( B^{pu})^T,
\qquad
E^{up}=(C^{pu})^T.
\end{aligned}
\]
The transpose exchanges the complete trial and test index groups,
including the continuum, component, and gradient-direction indices,
whenever applicable.

All effective coefficients are assumed to be independent of time.

\item[\textbf{(A2)}]
The tensors $A^{uu}(x)$, $A^{pp}(x)$, $H^{pp}(x)$,
$C^{uu}(x)$, $C^{pp}(x)$, and $D^{pp}(x)$ are uniformly positive definite a.e. in
$\Omega$. In particular, there exist constants
$\alpha_{uu},\alpha_{pp},\alpha_H,\alpha_{cu},\alpha_{cp},\alpha_{D}>0$ such that
\begin{align*}
\xi^T A^{uu}(x)\xi
&\ge \alpha_{uu}|\xi|^2,
&&\forall\,\xi\in\mathbb R^{d^2N_c},
\\
\eta^T A^{pp}(x)\eta
&\ge \alpha_{pp}|\eta|^2,
&&\forall\,\eta\in\mathbb R^{dN_c},
\\
\zeta^T H^{pp}(x)\zeta
&\ge \alpha_H|\zeta|^2,
&&\forall\,\zeta\in\mathbb R^{N_c},
\\
v^T C^{uu}(x)v
&\ge \alpha_{cu}|v|^2,
&&\forall\,v\in\mathbb R^{dN_c},
\\
\zeta^T C^{pp}(x)\zeta
&\ge \alpha_{cp}|\zeta|^2,
&&\forall\,\zeta\in\mathbb R^{N_c},
\\
\eta^T D^{pp}(x)\eta
&\ge \alpha_D|\eta|^2,
&&\forall\,\eta\in\mathbb R^{dN_c}.
\end{align*}

\item[\textbf{(A3)}]
There exist constants $\beta_{uu},\beta_{pp}, \beta_{gp} \ge 0$ such that,
for all $U\in V_u$ and $P\in V_p$,
\begin{align*}
2\left|(B^{uu}\nabla U,U)\right|
&\le \beta_{uu}\left(\|\nabla U\|^2+\|U\|^2\right),\\
2\left|(B^{pp}\nabla P,P)\right|
&\le \beta_{pp}\left(\|\nabla P\|^2+\|P\|^2\right),\\
2\left|(G^{pp}\nabla P,P)\right|&\le \beta_{gp}\left(\| \nabla P\|^2+\|P\|^2\right),
\end{align*}
where
\[
\beta_{uu}<\min\{\alpha_{uu},\alpha_{cu}\},
\qquad
\beta_{pp}<\min\{\alpha_{pp},\alpha_{cp}\},\qquad
\beta_{gp}<\min\{\alpha_{H},\alpha_{D}\}.
\]

\item[\textbf{(A4)}]
The source terms satisfy
\[
F_u\in L^{\infty}(0,T;L^2)^{dN_c},
\qquad
\partial_t F_u\in L^{2}(0,T;L^2)^{dN_c},
\qquad
F_p\in L^2(0,T;L^2),
\qquad
P_0\in V_p.
\]
In particular, these conditions imply
\[
F_u\in H^1(0,T;V_u'),
\qquad
F_p\in L^2(0,T;V_p').
\]

\item[\textbf{(A5)}]
We assume that the exact solution satisfies
\[
\begin{aligned}
U&\in
H^1(0,T;H^{k+1}(\Omega)^{dN_c})
\cap H^2(0,T;V_u),
\\
P&\in
H^1(0,T;H^{k+1}(\Omega)^{N_c})
\cap H^2(0,T;V_p).
\end{aligned}
\]
\end{itemize}
\end{assumption}

\begin{remark}
For the simplified model \eqref{eq:multicontinum}, we set
\[
A^{pu}=A^{up}=0,
\qquad
G^{pp}=\bar G^{pp}=0,
\qquad
D^{pp}=0.
\]
The conditions on these coefficients in \textbf{(A1)}
are then automatic. The positive definiteness requirement
for $D^{pp}$ in \textbf{(A2)}, and the estimate involving
$G^{pp}$ together with the restriction on $\beta_{gp}$
in \textbf{(A3)}, are removed.
All remaining assumptions are retained.
\end{remark}

\subsection{Weak formulation}
In this subsection, we present the weak formulation of the multicontinuum model \eqref{eq:multicontinum_full}.

Define
\[
\begin{aligned}
&V_u=[H_0^1(\Omega)]^{dN_c},
\qquad
V_p=[H_0^1(\Omega)]^{N_c},
\qquad
V=V_u\times V_p,\\
&\|U\|_{V_u}^2=\|\nabla U\|_{L^2}^2+\|U\|_{L^2}^2,
\qquad
\|P\|_{V_p}^2=\|\nabla P\|_{L^2}^2+\|P\|_{L^2}^2.
\end{aligned}
\]
For subsequent analysis, define the bilinear forms associated with the full
model \eqref{eq:multicontinum_full} by
\begin{align*}
\mathcal A_u(U,V)
={}&
(A^{uu}\nabla U,\nabla V)
+(B^{uu}\nabla U,V)
+(\bar B^{uu}U,\nabla V)
+(C^{uu}U,V),
\\
\mathcal A_p(P,Q)
={}&
(A^{pp}\nabla P,\nabla Q)
+(B^{pp}\nabla P,Q)
+(\bar B^{pp}P,\nabla Q)
+(C^{pp}P,Q),
\\
\mathcal B_{up}(P,V)
={}&
(A^{pu}\nabla P,\nabla V)
+(B^{pu}\nabla P,V)
+(\bar B^{pu}P,\nabla V)
+(C^{pu}P,V),
\\
\mathcal B_{pu}(Q,U)
={}&
(A^{up}\nabla U,\nabla Q)
+(G^{up}\nabla U,Q)
+(\bar G^{up}U,\nabla Q)
+(E^{up}U,Q),
\\
\mathcal M_p(P,Q)
={}&
(H^{pp}P,Q)
+(G^{pp}\nabla P,Q)
+(\bar G^{pp}P,\nabla Q)
+(D^{pp}\nabla P,\nabla Q).
\end{align*}

By the symmetry and reciprocity relations in Assumption \ref{assump:coeff_reg}, we have
\[
\mathcal B_{up}(P,V)
=
\mathcal B_{pu}(V,P),
\qquad P\in V_p,\quad V\in V_u.
\]
Therefore, for notational convenience, define
$\mathcal B(P,V)
:=
\mathcal B_{up}(P,V)
=
\mathcal B_{pu}(P,V).$

The weak formulation of \eqref{eq:multicontinum_full} is: find
$(U,P)\in V$ with the required time regularity such that, for all
$(V,Q)\in V$ and almost every $t\in(0,T)$,
\begin{equation}
\begin{aligned}
\mathcal A_u(U,V)-\mathcal B(P,V)
&=( F_u,V),
\\
\mathcal M_p(\partial_t P,Q)
+\mathcal B(Q,\partial_t U)
+\mathcal A_p(P,Q)
&=(F_p,Q).
\end{aligned}
\label{eq:weak_full_compact}
\end{equation}
The initial displacement $U_0\in V_u$ is defined by
\begin{equation}
\mathcal A_u(U_0,V)
-
\mathcal B(P_0,V)
=
(F_u(0),V),
\qquad V\in V_u.
\label{eq:continuous_initial_equilibrium_full}
\end{equation}

Equivalently, in fully expanded form,
\begin{equation*}
\begin{aligned}
&
(A^{uu}\nabla U,\nabla V)
+(B^{uu}\nabla U,V)
+(\bar B^{uu}U,\nabla V)
+(C^{uu}U,V)
\\
&\qquad
-(A^{pu}\nabla P,\nabla V)
-(B^{pu}\nabla P,V)
-(\bar B^{pu}P,\nabla V)
-(C^{pu}P,V)
=
(F_u,V),
\\[0.4em]
&
(H^{pp}\partial_t P,Q)
+(G^{pp}\nabla\partial_t P,Q)
+(\bar G^{pp}\partial_t P,\nabla Q)
+(D^{pp}\nabla\partial_t P,\nabla Q)
\\
&\qquad
+(A^{up}\nabla\partial_t U,\nabla Q)
+(G^{up}\nabla\partial_t U,Q)
+(\bar G^{up}\partial_t U,\nabla Q)
+(E^{up}\partial_t U,Q)
\\
&\qquad
+(A^{pp}\nabla P,\nabla Q)
+(B^{pp}\nabla P,Q)
+(\bar B^{pp}P,\nabla Q)
+(C^{pp}P,Q)
=
(F_p,Q).
\end{aligned}
\label{eq:macro_weak_full_expanded}
\end{equation*}

\subsection{Fully discrete scheme}
We assume a uniform time partition $t_n=n\tau,~ n=0,\ldots,N,~ \tau=T/N$. 
For any function $z$, set $D_tz^{n+1}=\frac{z^{n+1}-z^n}{\tau}$, where $z^n=z(t_n)$.

Let
$V_{u,h}\subset V_u$ and $V_{p,h}\subset V_p$ be conforming degree-$k$
Lagrange finite element spaces, with $k\ge1$.

For the fully discrete scheme, we use the backward Euler method for
time discretization.
Given $(U_h^n,P_h^n)$, find
$(U_h^{n+1},P_h^{n+1})\in V_{u,h}\times V_{p,h}$ such that, for all
$(V_h,Q_h)\in V_{u,h}\times V_{p,h}$,
\begin{equation}
\begin{aligned}
\mathcal A_u(U_h^{n+1},V_h)
-
\mathcal B(P_h^{n+1},V_h)
&=
( F_u^{n+1},V_h),
\\
\mathcal M_p(D_tP_h^{n+1},Q_h)
+
\mathcal B(Q_h,D_tU_h^{n+1})
+
\mathcal A_p(P_h^{n+1},Q_h)
&=
(F_p^{n+1},Q_h) .
\end{aligned}
\label{eq:fully_discrete_compact_full}
\end{equation}
Equivalently, in fully expanded form,
\begin{equation*}
\begin{aligned}
&
(A^{uu}\nabla U_h^{n+1},\nabla V_h)
+
(B^{uu}\nabla U_h^{n+1},V_h)
+
(\bar B^{uu}U_h^{n+1},\nabla V_h)
+
(C^{uu}U_h^{n+1},V_h)
\\
&\qquad
-
(A^{pu}\nabla P_h^{n+1},\nabla V_h)
-
(B^{pu}\nabla P_h^{n+1},V_h)
-
(\bar B^{pu}P_h^{n+1},\nabla V_h)
-
(C^{pu}P_h^{n+1},V_h)
\\
&\qquad
=
(F_u^{n+1},V_h),
\\[0.6em]
&
(H^{pp}D_tP_h^{n+1},Q_h)
+
(G^{pp}\nabla D_tP_h^{n+1},Q_h)
+
(\bar G^{pp}D_tP_h^{n+1},\nabla Q_h)
+
(D^{pp}\nabla D_tP_h^{n+1},\nabla Q_h)
\\
&\qquad
+
(A^{up}\nabla D_tU_h^{n+1},\nabla Q_h)
+
(G^{up}\nabla D_tU_h^{n+1},Q_h)
+
(\bar G^{up}D_tU_h^{n+1},\nabla Q_h)
+
(E^{up}D_tU_h^{n+1},Q_h)
\\
&\qquad
+
(A^{pp}\nabla P_h^{n+1},\nabla Q_h)
+
(B^{pp}\nabla P_h^{n+1},Q_h)
+
(\bar B^{pp}P_h^{n+1},\nabla Q_h)
+
(C^{pp}P_h^{n+1},Q_h)
\\
&\qquad
=
(F_p^{n+1},Q_h) .
\end{aligned}
\label{eq:fully_discrete_full}
\end{equation*}

The initial pressure is taken as
\[
P_h^0=\Pi_hP_0,
\]
where $\Pi_h$ is the $L^2$ projection onto $V_{p,h}$, and $U_h^0$ is
defined by
\begin{equation}
\begin{aligned}
\mathcal A_u(U_h^{0},V_h)
-
\mathcal B(P_h^{0},V_h)
=
(F_u^{0},V_h),\qquad V_h\in V_{u,h}.
\end{aligned}
\label{eq:macro_initial_displacement}
\end{equation}


\section{Well-posedness and stability}
\label{sec:wellposedness}

In this section, we establish the well-posedness and stability of the weak formulation
\eqref{eq:weak_full_compact} and the fully discrete scheme
\eqref{eq:fully_discrete_compact_full}, separately.

\begin{lemma}
\label{lem:diagonal_coercivity_full}
Suppose Assumption \ref{assump:coeff_reg} \textbf{(A1)-(A3)} hold.
There exist constants
$c_u,c_p>0$, independent of the discretization parameters, such that
\begin{align}
\mathcal A_u(U,U)
&\ge c_u\|U\|_{V_u}^2,
&&U\in V_u,
\label{eq:Au_coercive_full}
\\
\mathcal A_p(P,P)
&\ge c_p\|P\|_{V_p}^2,
&&P\in V_p,
\label{eq:Ap_coercive_full}
\\
\mathcal M_p(P,P)
&\ge c_m\|P\|_{V_p}^2,
&&P\in V_p.
\label{eq:Mp_coercive_full}
\end{align}
Moreover, 
\begin{align}
\mathcal A_u(U_h^{n+1},D_t U_h^{n+1})
&\ge \frac{c_u}{2} D_t\|U_h^{n+1}\|_{V_u}^2,
&&U_h^{n+1}\in V_{u,h},
\label{eq:Au_coercive_full_discrete}
\\
\mathcal M_p(D_t P_h^{n+1},P_h^{n+1})
&\ge \frac{c_m}{2}D_t \|P_h^{n+1}\|_{V_p}^2,
&&P_h^{n+1}\in V_{p,h},
\label{eq:Mp_coercive_full_discrete}
\end{align}
where
\[
\begin{aligned}
c_u&=\min\{\alpha_{uu},\alpha_{cu}\}-\beta_{uu}>0,
\\
c_p&=\min\{\alpha_{pp},\alpha_{cp}\}-\beta_{pp}>0,
\\
c_m&=\min\{\alpha_H,\alpha_D\}-\beta_{gp}>0.
\end{aligned}
\]
\end{lemma}
\begin{proof}
By the symmetry relations in \textbf{(A1)},
\[
(\bar B^{uu}U,\nabla U)=(B^{uu}\nabla U,U),
\qquad
(\bar B^{pp}P,\nabla P)=(B^{pp}\nabla P,P),
\qquad
(\bar G^{pp}P,\nabla P)=(G^{pp}\nabla P,P).
\]
By the definition of $\mathcal A_u$ and the relation
$\bar B^{uu}=(B^{uu})^T$, we have
\[
\begin{aligned}
\mathcal A_u(U,U)
={}&(A^{uu}\nabla U,\nabla U)
+2(B^{uu}\nabla U,U)+(C^{uu}U,U).
\end{aligned}
\]
Consequently, \textbf{(A2)}-\textbf{(A3)} imply
\[
\begin{aligned}
&(A^{uu}\nabla U,\nabla U)
\ge \alpha_{uu}\|\nabla U\|^2,\quad
(C^{uu}U,U)
\ge \alpha_{cu}\|U\|^2,\\
&|2(B^{uu}\nabla U,U)|
\ge -b_{uu}\|\nabla U\|^2-b_{uu}\|U\|^2.
\end{aligned}
\]
Combining these estimates gives
\[
\begin{aligned}
\mathcal A_u(U,U)
&\ge (\alpha_{uu}-b_{uu})\|\nabla U\|^2+(\alpha_{cu}-b_{uu})\|U\|^2\\
&\ge c_u\left(\|\nabla U\|^2+\|U\|^2\right) \ge c_u\|U\|_{V_u}^2.
\end{aligned}
\]
where
$c_u=\min\{\alpha_{uu},\alpha_{cu}\}-\beta_{uu}>0$.

Besides, we have
\[
\begin{aligned}
&\mathcal A_u(U_h^{n+1},D_t U_h^{n+1})\\
={}&(A^{uu}\nabla U_h^{n+1},\nabla D_tU_h^{n+1})
+(B^{uu}\nabla U_h^{n+1},D_t U_h^{n+1})+(\bar B^{uu}U_h^{n+1},\nabla D_t U_h^{n+1})+(C^{uu}U_h^{n+1},D_t U_h^{n+1})\\
= &\frac12D_t\left[(A^{uu}\nabla U_h^{n+1},\nabla U_h^{n+1})+(C^{uu}U_h^{n+1},U_h^{n+1})\right]
+ D_t(B^{uu}\nabla U_h^{n+1},U_h^{n+1})\\
\ge & \frac{c_u}{2} D_t\|U_h^{n+1}\|_{V_u}^2.  
\end{aligned}
\]

The proof for $\mathcal A_p$, $\mathcal M_p$ are similar. This completes the proof.
\end{proof}

Next, we establish the existence and uniqueness of the weak solution to the full multicontinuum model \eqref{eq:multicontinum_full}.

\begin{theorem}
\label{thm:weak_existence_full}
Suppose Assumption \ref{assump:coeff_reg} \textbf{(A1)-(A4)} hold.
Then the weak problem \eqref{eq:weak_full_compact} has a unique solution $(U,P)$ such that
\[
U\in H^1(0,T;V_u),
\qquad
P\in H^1(0,T;V_p).
\]
\end{theorem}

\begin{proof}
Define the linear operator $\mathcal T:V_p\to V_u$ by
\begin{equation}
\mathcal A_u(\mathcal TR,V)=\mathcal B(R,V),
\qquad V\in V_u.
\label{eq:response_operator_full}
\end{equation}
By the continuity and coercivity of $\mathcal A_u$, and $|\mathcal B(R,V)|
\leq C_B\|R\|_{V_p}\|V\|_{V_u}$, the Lax--Milgram theorem implies that $\mathcal T$ is
well defined and bounded. In particular,
\[
\|\mathcal TR\|_{V_u}\leq C\|R\|_{V_p}.
\]

Next, for almost every $t\in(0,T)$, let $w(t)\in V_u$ solve
\begin{equation}
\mathcal A_u(w(t),V)= (F_u(t),V),
\qquad V\in V_u.
\label{eq:load_response_full}
\end{equation}
The Lax-Milgram theorem gives a unique $w(t)$. 
By Assumptions \textbf{(A4)}, we have $F_u\in H^1(0,T;V_u')$, and 
since the effective
coefficients are independent of time, the solution operator associated with
$\mathcal A_u$ is time independent and bounded, consequently,
\[
w\in H^1(0,T;V_u).
\]
The first equation of \eqref{eq:weak_full_compact} therefore yields
\[
\mathcal A_u(U-w,V)=\mathcal B(P,V)
\qquad\forall V\in V_u,
\]
which, by the definition of $\mathcal T$, is equivalent to
\begin{equation}
U=w+\mathcal TP.
\label{eq:displacement_elimination_full}
\end{equation}
In particular,
$\partial_tU=\partial_tw+\mathcal T\partial_tP$.

For $R,Q\in V_p$, define
\begin{equation}
\mathcal S(R,Q)
:=\mathcal M_p(R,Q)
+\mathcal A_u(\mathcal TR,\mathcal TQ).
\label{eq:reduced_storage_form_full}
\end{equation}
By Assumption \textbf{(A1)}, the forms $\mathcal A_u$ and
$\mathcal M_p$ are symmetric. Hence $\mathcal S$ is symmetric.
Moreover, the continuity of $\mathcal M_p$, $\mathcal A_u$, and
$\mathcal T$ implies that $\mathcal S$ is continuous on
$V_p\times V_p$. By the coercivity of $\mathcal M_p$,
\[
\begin{aligned}
\mathcal S(Q,Q)
&=\mathcal M_p(Q,Q)
  +\mathcal A_u(\mathcal TQ,\mathcal TQ)\\
&\geq \mathcal M_p(Q,Q)
 \geq c_m\|Q\|_{V_p}^2.
\end{aligned}
\]
Thus, $\mathcal S$ is coercive on $V_p$ and defines an equivalent
inner product on this space.

By the definition of $\mathcal T$ and the symmetry of $\mathcal A_u$,
\[
\begin{aligned}
\mathcal B(Q,\mathcal T\partial_tP)
=\mathcal A_u(\mathcal TQ,\mathcal T\partial_tP)
=\mathcal A_u(\mathcal T\partial_tP,\mathcal TQ).
\end{aligned}
\]
Substituting \eqref{eq:displacement_elimination_full} into the second
equation of \eqref{eq:weak_full_compact} gives the reduced pressure
equation
\begin{equation}
\mathcal S(\partial_tP,Q)+\mathcal A_p(P,Q)
= (F_p,Q) -\mathcal B(Q,\partial_tw),
\qquad Q\in V_p.
\label{eq:reduced_pressure_equation_full}
\end{equation}

Let $\mathbf S:V_p\to V_p'$ and
$\mathbf A_p:V_p\to V_p'$ be the operators induced by
$\mathcal S$ and $\mathcal A_p$, respectively; that is,
\[
\langle\mathbf SR,Q\rangle=\mathcal S(R,Q),
\qquad
\langle\mathbf A_pR,Q\rangle=\mathcal A_p(R,Q).
\]
The continuity and coercivity of $\mathcal S$ imply, by the
Lax-Milgram theorem, that $\mathbf S$ is an isomorphism from
$V_p$ onto $V_p'$.

Define $\ell\in L^2(0,T;V_p')$ by
\[
\langle\ell(t),Q\rangle
:= (F_p(t),Q)
  -\mathcal B(Q,\partial_tw(t)),
\qquad Q\in V_p.
\]
Indeed, the continuity of $\mathcal B$, together with
$F_p\in L^2(0,T;V_p')$ and
$\partial_tw\in L^2(0,T;V_u)$, gives
$\ell\in L^2(0,T;V_p')$.

In operator form, \eqref{eq:reduced_pressure_equation_full} becomes
\begin{equation}
\partial_tP+\mathbf S^{-1}\mathbf A_pP
=\mathbf S^{-1}\ell,
\qquad P(0)=P_0.
\label{eq:abstract_reduced_evolution}
\end{equation}
Since $\mathbf A_p:V_p\to V_p'$ is bounded and
$\mathbf S^{-1}:V_p'\to V_p$ is bounded,
$\mathbf S^{-1}\mathbf A_p$ is a bounded linear operator on $V_p$.
Furthermore, $\mathbf S^{-1}\ell\in L^2(0,T;V_p)$.
The standard existence and uniqueness result for linear evolution
equations with bounded operators (equivalently, the variation-of-constants
formula) therefore yields a unique solution
\[
P\in H^1(0,T;V_p)
\]
with the prescribed initial value $P(0)=P_0$.

Define $U$ by \eqref{eq:displacement_elimination_full}. Since
$w\in H^1(0,T;V_u)$, $\mathcal T$ is bounded, and
$P\in H^1(0,T;V_p)$, it follows that
\[
U\in H^1(0,T;V_u).
\]
The construction of $P$ and $U$ shows that both equations in
\eqref{eq:weak_full_compact} are satisfied for almost every
$t\in(0,T)$.

It remains to verify the initial equilibrium condition. Since
$w\in H^1(0,T;V_u)$ and $P\in H^1(0,T;V_p)$, both have well-defined
traces at $t=0$. Using $P(0)=P_0$, the definition of $w$, and that of
$\mathcal T$, we obtain, for every $V\in V_u$,
\[
\begin{aligned}
\mathcal A_u(U(0),V)-\mathcal B(P_0,V)
&=\mathcal A_u(w(0)+\mathcal TP_0,V)
  -\mathcal B(P_0,V)\\
&=(F_u(0),V)
  +\mathcal B(P_0,V)-\mathcal B(P_0,V)\\
&=(F_u(0),V).
\end{aligned}
\]
Thus, the initial equilibrium condition
\eqref{eq:continuous_initial_equilibrium_full} holds.

Finally, the reduced evolution equation \eqref{eq:abstract_reduced_evolution}
has a unique solution for the prescribed initial value. Hence the pressure
$P$ is unique, and the representation
$U=w+\mathcal TP$ then implies the uniqueness of $U$.
This completes the proof.
\end{proof}

Then we give the energy identity and stability of the weak solution to the full multicontinuum model \eqref{eq:multicontinum_full}.

\begin{theorem}
\label{thm:weak_stability_full}
Suppose Assumption \ref{assump:coeff_reg} \textbf{(A1)-(A4)} hold.
The solution of the weak problem \eqref{eq:weak_full_compact} satisfies
\begin{equation}
\begin{aligned}
&
\sup_{0\le t\le T}
\left(
\|U(t)\|_{V_u}^2+\|P(t)\|_{V_p}^2
\right)
+
\int_0^T\|P(t)\|_{V_p}^2\,dt
\\
&\qquad\le
C\left(
\|U_0\|_{V_u}^2+\|P_0\|_{V_p}^2
+
\|F_u\|_{L^{\infty}(0,T;L^2)}^2
+
\|\partial _t F_u\|_{L^2(0,T;L^2)}^2
+
\|F_p\|_{L^2(0,T;L^2)}^2
\right).
\end{aligned}
\label{eq:weak_stability_bound_full}
\end{equation}
\end{theorem}

\begin{proof}
Take $V=\partial_tU$ and $Q=P$ in
\eqref{eq:weak_full_compact}.
The coupling contributions cancel:
\[
-\mathcal B(P,\partial_tU)
+
\mathcal B(P,\partial_tU)=0.
\]
Since the forms $\mathcal A_u$ and $\mathcal M_p$ are
symmetric and independent of time, we get
\[
\mathcal A_u(U,\partial_tU)
=
\frac12\frac{d}{dt}\mathcal A_u(U,U),
\]
\[
\mathcal M_p(\partial_tP,P)
=
\frac12\frac{d}{dt}\mathcal M_p(P,P).
\]
Define the energy by
\begin{equation}
\mathcal E(t)
=
\frac12\mathcal A_u(U(t),U(t))
+
\frac12\mathcal M_p(P(t),P(t)).
\label{eq:continuous_energy_full}
\end{equation}
Then we have
\begin{equation}
\frac{d}{dt}\mathcal E(t)
+
\mathcal A_p(P(t),P(t))
=
( F_u(t),\partial_tU(t))
+
( F_p(t),P(t))
\label{eq:energy_identity_continuous_full}
\end{equation}

Integrating in time and integrating the displacement load
by parts yields
\begin{equation}
\begin{aligned}
\mathcal E(t)
+
\int_0^t\mathcal A_p(P,P)\,ds
={}&
\mathcal E(0)
+
( F_u(t),U(t))
-
 (F_u(0),U_0)
\\
&-
\int_0^t (\partial_tF_u,U) \,ds
+
\int_0^t ( F_p,P) \,ds.
\end{aligned}
\label{eq:integrated_energy_full}
\end{equation}
By Cauchy-Schwarz's inequality and Young's inequality,
\[
\begin{aligned}
|  (F_u(t),U(t)) |
&\le
\varepsilon \|U(t)\|_{L^2}^2
+
C_\varepsilon\|F_u(t)\|_{L^2}^2,\\
|  (F_u(0),U_0) |
&\le
C\left(
\|F_u(0)\|_{L^2}^2+\|U_0\|_{L^2}^2
\right),\\
\int_0^t|( \partial_tF_u,U) |\,ds
&\le
\varepsilon \|U\|_{L^2(0,T;L^2)}^2
+
C_\varepsilon 
\|\partial_tF_u\|_{L^2(0,T;L^2)}^2,\\
\int_0^t| ( F_p,P) |\,ds
&\le
\varepsilon\|P\|_{L^2(0,T;L^2)}^2
+
C_\varepsilon \|F_p\|_{L^2(0,T;L^2)}^2.
\end{aligned}
\]
Continuity and coercivity imply
\begin{equation*}
c_E\bigl(\|U\|_{V_u}^2+\|P\|_{V_p}^2\bigr)
\le \mathcal E
\le
C_E\bigl(\|U\|_{V_u}^2+\|P\|_{V_p}^2\bigr),
\label{eq:energy_equivalence_full}
\end{equation*}
where $c_E,C_E>0$.

Insert these estimates into
\eqref{eq:integrated_energy_full}, take the supremum over
$t$, and choose $\varepsilon$ sufficiently small to absorb
the terms involving $\sup_{0\le t\le T}\mathcal E(t)$.
Returning to \eqref{eq:integrated_energy_full} at $t=T$
then bounds the dissipation integral. The proof is complete.
\end{proof}

Finally, we establish the well-posedness and stability of the fully discrete scheme \eqref{eq:fully_discrete_compact_full}.

\begin{theorem}
\label{thm:fd_wellposed}
Suppose Assumption \ref{assump:coeff_reg} \textbf{(A1)-(A4)} hold. The fully discrete scheme
\eqref{eq:fully_discrete_compact_full} has a unique solution
$(U_h^{n+1},P_h^{n+1})\in V_{u,h}\times V_{p,h}$. Moreover, the solution
satisfies a discrete stability estimate of the form
\begin{equation}
\begin{aligned}
&
\max_{0\le n\le N}
\left(
\|U_h^n\|_{V_u}^2+\|P_h^n\|_{V_p}^2
\right)+
\tau\sum_{n=0}^{N-1}\|P_h^{n+1}\|_{V_p}^2
\\
\le &
C\Bigg[
\|U_h^0\|_{V_u}^2+\|P_h^0\|_{V_p}^2
+
\|F_u\|_{L^{\infty}(0,T;L^2)}^2+\|\partial_tF_u\|_{L^2(0,T;L^2)}^2+\|F_p\|_{L^2(0,T;L^2)}^2
\Bigg],
\end{aligned}
\label{eq:discrete_stability_full}
\end{equation}
where $C$ is independent of $h$ and $\tau$.
\end{theorem}
\begin{proof}
First, we prove the well-poseness of the fully discrete scheme
\eqref{eq:fully_discrete_compact_full}.

Consider the difference $(\widehat U,\widehat P)$ of two
solutions at the same time step, with identical previous
values and loads.
Then
\[
\mathcal A_u(\widehat U,V_h)
-
\mathcal B(\widehat P,V_h)=0,
\]
after multiplying the pressure equation by $\tau$,
\[
\mathcal M_p(\widehat P,Q_h)
+
\mathcal B(Q_h,\widehat U)
+
\tau\mathcal A_p(\widehat P,Q_h)=0.
\]
Take $V_h=\widehat U$ and $Q_h=\widehat P$ and add the two equations, we get
\[
\mathcal A_u(\widehat U,\widehat U)
+
\mathcal M_p(\widehat P,\widehat P)
+
\tau\mathcal A_p(\widehat P,\widehat P)=0,
\]
which implies $\widehat U=0$ and $\widehat P=0$.
Since the algebraic system is square and finite dimensional,
uniqueness implies existence.

Next, we prove the energy estimate of the fully discrete scheme
\eqref{eq:fully_discrete_compact_full}.

Take
$V_h=D_t U_h^{n+1},
~
Q_h=P_h^{n+1}$ in \eqref{eq:fully_discrete_compact_full},
and add the two
equations.
The coupling terms cancel exactly. 

By Lemma \ref{lem:diagonal_coercivity_full}, we have
\[
\mathcal A_u(U_h^{n+1},D_t U_h^{n+1})
\ge \frac{c_u}{2} D_t\|U_h^{n+1}\|_{V_u}^2,
\qquad
\mathcal M_p(D_t P_h^{n+1},P_h^{n+1})
\ge \frac{c_m}{2}\|P_h^{n+1}\|_{V_p}^2,
\]
Then, the system \eqref{eq:fully_discrete_compact_full} gives
\begin{equation}
\begin{aligned}
D_t\|U_h^{n+1}\|_{V_u}^2+D_t \|P_h^{n+1}\|_{V_p}^2
+\|P_h^{n+1}\|_{V_p}^2
\le
( F_u^{n+1},D_t U_h^{n+1})
+
(F_p^{n+1},P_h^{n+1}).
\end{aligned}
\label{eq:discrete_energy_identity_full}
\end{equation}
Sum \eqref{eq:discrete_energy_identity_full} from $n=0$
to $N-1$, we get
\begin{equation}
\begin{aligned}
\|U_h^{N}\|_{V_u}^2+ \|P_h^{N}\|_{V_p}^2
+\sum_{n=0}^{N-1}\|P_h^{n+1}\|_{V_p}^2
\le
\sum_{n=0}^{N-1}\tau( F_u^{n+1},D_t U_h^{n+1})
+
\sum_{n=0}^{N-1}\tau(F_p^{n+1},P_h^{n+1}).
\end{aligned}
\label{eq:discrete_energy_identity_full2}
\end{equation}
By Cauchy-Schwarz and Young's inequalities,
\[
\begin{aligned}
\sum_{n=0}^{N-1}\tau
( F_u^{n+1},D_t U_h^{n+1})
={}&
( F_u^N,U_h^N)
-
( F_u^0,U_h^0)
-
\sum_{n=0}^{N-1}\tau
( D_tF_u^{n+1},U_h^n)\\
\le & \varepsilon
\max_{0\le n\le N}\|U_h^n\|_{V_u}^2
+
C_\varepsilon\left[
\|U_h^0\|_{V_u}^2
+
\max_{0\le n\le N}\|F_u^n\|_{L^2}^2
+
\sum_{n=0}^{N-1}\tau
\|D_tF_u^{n+1}\|_{L^2}^2
\right]\\
\le & \varepsilon
\max_{0\le n\le N}\|U_h^n\|_{L^2}^2
+
C_\varepsilon\left[
\|U_h^0\|_{L^2}^2
+
\max_{0\le n\le N}\|F_u^n\|_{L^2}^2
+
\int_{0}^{T}\|\partial_tF_u\|_{L^2}^2dt
\right].
\end{aligned}
\]
For the pressure,
\[
\sum_{n=0}^{N-1}\tau| (F_p^{n+1},P_h^{n+1})|
\le
\varepsilon \sum_{n=0}^{N-1}\tau \|P_h^{n+1}\|_{V_p}^2
+
C_\varepsilon\tau\|F_p^{n+1}\|_{L^2}^2.
\]
By choosing $\varepsilon$ sufficiently small, we can absorb the terms involving $\|U_h^n\|_{V_u}^2$ and $\|P_h^{n+1}\|_{V_p}^2$ into the left-hand side of the inequality \eqref{eq:discrete_energy_identity_full2}, yielding the desired stability estimate \eqref{eq:discrete_stability_full}.
\end{proof}

\section{Error analysis for the fully discrete full model}
\label{sec:error}

\subsection{Finite element projections}
In this subsection, we define the finite element projections and establish their approximation properties.

For $U\in V_u$, define $I_hU\in V_{u,h}$ by
\begin{equation}
\mathcal A_u(I_hU,V_h)
=
\mathcal A_u(U,V_h),
\qquad V_h\in V_{u,h}.
\label{eq:elliptic_projection_U_full}
\end{equation}
For $P\in H_p$, define $\Pi_hP\in V_{p,h}$ by
\begin{equation}
(\Pi_hP-P,Q_h)=0,
\qquad Q_h\in V_{p,h}.
\label{eq:L2_projection_P_full}
\end{equation}

\begin{lemma}[Projection estimates \cite{wheeler1973priori}]
\label{lem:projection_full}
For
\[
U\in H^{k+1}(\Omega)^{dN_c}\cap V_u,
\qquad
P\in H^{k+1}(\Omega)^{N_c}\cap V_p,
\]
the projections satisfy
\begin{equation}
\|U-I_hU\|_{V_u}
\le Ch^k\|U\|_{H^{k+1}},
\label{eq:U_projection_bound_full}
\end{equation}
and
\begin{equation}
\|P-\Pi_hP\|_{L^2}
+
h\|P-\Pi_hP\|_{V_p}
\le Ch^{k+1}\|P\|_{H^{k+1}}.
\label{eq:P_projection_bound_full}
\end{equation}
\end{lemma}

\begin{proof}
The estimate for $I_h$ follows from the continuity and
coercivity of $\mathcal A_u$, Cea's lemma, and standard
finite element approximation.

The $L^2$ estimate for $\Pi_h$ follows from its
best-approximation property.
On a shape-regular, quasi-uniform mesh family, the
$L^2$ projection is uniformly stable in $H^1$.
Combining this stability with a standard conforming
interpolant gives the $H^1$ estimate.
\end{proof}

\subsection{A priori error estimate}

\begin{theorem}
\label{thm:fd_error_full}
Suppose Assumption \ref{assump:coeff_reg} holds.
Let $(U_h^n,P_h^n)$ be the solution of
\eqref{eq:fully_discrete_compact_full}, $(U,P)$ be the solution of
\eqref{eq:weak_full_compact}.
Then
\begin{equation}
\begin{aligned}
\max_{0\le n\le N}
\left(
\|U^n-U_h^n\|_{V_u}^2+\|P^n-P_h^n\|_{V_p}^2
\right)
+
\tau\sum_{n=0}^{N-1}\|e_p^{n+1}\|_{V_p}^2
\le
C(
h^{2k}
+
\tau^2),
\end{aligned}
\label{eq:error_estimate_full}
\end{equation}
where the constant $C$ is independent of $h$ and $\tau$.
\end{theorem}

\begin{proof}
Set
\[
\begin{aligned}
e_u^n&=\eta_u^n+\theta_u^n,
&
e_p^n&=\eta_p^n+\theta_p^n,\\
\eta_u^n&=U(t_n)-I_hU(t_n),
&
\theta_u^n&=I_hU(t_n)-U_h^n,
\\
\eta_p^n&=P(t_n)-\Pi_hP(t_n),
&
\theta_p^n&=\Pi_hP(t_n)-P_h^n.
\end{aligned}
\]

Define
\[
\begin{aligned}
&\rho_u^{n+1}
=
\partial_tU(t_{n+1})-D_tU(t_{n+1}),\\
&\rho_p^{n+1}
=
\partial_tP(t_{n+1})-D_tP(t_{n+1}).
\end{aligned}
\]
Consequently,
\begin{equation}
\begin{aligned}
\tau\sum_{n=0}^{N-1}
\|\rho_u^{n+1}\|_{V_u}^2
&\le
C\tau^2
\|\partial_{tt}U\|_{L^2(0,T;V_u)}^2,
\\
\tau\sum_{n=0}^{N-1}
\|\rho_p^{n+1}\|_{V_p}^2
&\le
C\tau^2
\|\partial_{tt}P\|_{L^2(0,T;V_p)}^2.
\end{aligned}
\label{eq:time_consistency_bounds_full}
\end{equation}

Subtracting the discrete equations from \eqref{eq:fully_discrete_compact_full} the continuous
equations \eqref{eq:weak_full_compact} at $t_{n+1}$ gives
\begin{equation}
\mathcal A_u(e_u^{n+1},V_h)
-
\mathcal B(e_p^{n+1},V_h)=0,
\label{eq:total_error_u_full}
\end{equation}
\begin{equation}
\begin{aligned}
\mathcal M_p(D_te_p^{n+1},Q_h)
+
\mathcal B(Q_h,D_te_u^{n+1})
+
\mathcal A_p(e_p^{n+1},Q_h)
=
-\mathcal M_p(\rho_p^{n+1},Q_h)
-\mathcal B(Q_h,\rho_u^{n+1}).
\end{aligned}
\label{eq:total_error_p_full}
\end{equation}
Using the error decomposition in
\eqref{eq:total_error_u_full}-\eqref{eq:total_error_p_full}, and the projection projection $\mathcal A_u(\eta_u^n,V_h)=0$,
we obtain
\begin{equation}
\begin{aligned}
\mathcal A_u(\theta_u^{n+1},V_h)
-
\mathcal B(\theta_p^{n+1},V_h)
&=
(R_u^{n+1},V_h),
\\
\mathcal M_p(D_t\theta_p^{n+1},Q_h)
+
\mathcal B(Q_h,D_t\theta_u^{n+1})
+
\mathcal A_p(\theta_p^{n+1},Q_h)
&=
(R_p^{n+1},Q_h),
\end{aligned}
\label{eq:projected_error_system_full}
\end{equation}
where
\begin{equation*}
 (R_u^n,V_h)
=
\mathcal B(\eta_p^n,V_h),
\label{eq:displacement_residual_full}
\end{equation*}
\begin{equation*}
\begin{aligned}
(R_p^{n+1},Q_h)
={}&
-\mathcal M_p(D_t\eta_p^{n+1},Q_h)
-\mathcal B(Q_h,D_t\eta_u^{n+1})
-\mathcal A_p(\eta_p^{n+1},Q_h)
\\
&-
\mathcal M_p(\rho_p^{n+1},Q_h)
-\mathcal B(Q_h,\rho_u^{n+1}).
\end{aligned}
\label{eq:pressure_residual_full}
\end{equation*}

The error system \eqref{eq:projected_error_system_full}
has the same structure as the full discrete scheme.
Applying Theorem \ref{thm:fd_wellposed} with loads
$R_u^n$ and $R_p^{n+1}$ yields
\begin{equation}
\begin{aligned}
&
\max_{0\le n\le N}
\left(
\|\theta_u^n\|_{V_u}^2
+
\|\theta_p^n\|_{V_p}^2
\right)
+
\tau\sum_{n=0}^{N-1}
\|\theta_p^{n+1}\|_{V_p}^2
\\
&\qquad\le
C\Bigg[
\|\theta_u^0\|_{V_u}^2
+
\|\theta_p^0\|_{V_p}^2
+
\max_{0\le n\le N}\|R_u^n\|_{L^2}^2
+
\tau\sum_{n=0}^{N-1}
\|D_tR_u^{n+1}\|_{L^2}^2
+
\tau\sum_{n=0}^{N-1}
\|R_p^{n+1}\|_{L^2}^2
\Bigg].
\end{aligned}
\label{eq:error_stability_with_residuals_full}
\end{equation}

Next, we bound the right-term of \eqref{eq:error_stability_with_residuals_full}.
Continuity of the coupling form gives
\[
\|R_u^n\|_{L^2}
\le C\|\eta_p^n\|_{V_p},
\]
and, because the coefficients are independent of time,
\[
\|D_tR_u^{n+1}\|_{L^2}
\le C\|D_t\eta_p^{n+1}\|_{V_p}.
\]
Moreover,
\[
\begin{aligned}
\|R_p^{n+1}\|_{L^2}
\le C\Big(
&
\|D_t\eta_p^{n+1}\|_{V_p}
+
\|D_t\eta_u^{n+1}\|_{V_u}
+
\|\eta_p^{n+1}\|_{V_p}
+
\|\rho_p^{n+1}\|_{V_p}
+
\|\rho_u^{n+1}\|_{V_u}
\Big).
\end{aligned}
\]
Since the projections are independent of time,
\[
D_t\eta_u^{n+1}
=
\frac1\tau
\int_{t_n}^{t_{n+1}}
(I-I_h)\partial_tU(s)\,ds,
\]
\[
D_t\eta_p^{n+1}
=
\frac1\tau
\int_{t_n}^{t_{n+1}}
(I-\Pi_h)\partial_tP(s)\,ds.
\]
Hence
\begin{equation}
\begin{aligned}
\tau\sum_{n=0}^{N-1}
\|D_t\eta_u^{n+1}\|_{V_u}^2
&\le
Ch^{2k}
\|\partial_tU\|_{L^2(0,T;H^{k+1})}^2,
\\
\tau\sum_{n=0}^{N-1}
\|D_t\eta_p^{n+1}\|_{V_p}^2
&\le
Ch^{2k}
\|\partial_tP\|_{L^2(0,T;H^{k+1})}^2.
\end{aligned}
\label{eq:projection_time_derivative_bounds_full}
\end{equation}
By Lemma \ref{lem:projection_full},
\[
\|\eta_u^n\|_{V_u}
\le Ch^k\|U(t_n)\|_{H^{k+1}},
\qquad
\|\eta_p^n\|_{V_p}
\le Ch^k\|P(t_n)\|_{H^{k+1}}.
\]
Combining these bounds with
\eqref{eq:time_consistency_bounds_full} yields
\begin{equation}
\begin{aligned}
&
\max_n\|R_u^n\|_{V_u'}^2
+
\tau\sum_{n=0}^{N-1}
\|D_tR_u^{n+1}\|_{V_u'}^2
+
\tau\sum_{n=0}^{N-1}
\|R_p^{n+1}\|_{V_p'}^2
\\
&\qquad\le
Ch^{2k}\left(
\|U\|_{H^1(0,T;H^{k+1})}^2
+
\|P\|_{H^1(0,T;H^{k+1})}^2
\right)
\\
&\hspace{4em}
+
C\tau^2\left(
\|\partial_{tt}U\|_{L^2(0,T;V_u)}^2
+
\|\partial_{tt}P\|_{L^2(0,T;V_p)}^2
\right).
\end{aligned}
\label{eq:residual_estimate_full}
\end{equation}

Since $P_h^0=\Pi_hP_0$, we have
$\theta_p^0=0.$
Subtracting the continuous and discrete equilibrium
equations at $t=0$ and using the elliptic projection
property gives
\[
\mathcal A_u(\theta_u^0,V_h)
=
\mathcal B(\eta_p^0,V_h).
\]
Taking $V_h=\theta_u^0$ and using coercivity yields
\[
\|\theta_u^0\|_{V_u}
\le
C\|\eta_p^0\|_{V_p}
\le
Ch^k\|P_0\|_{H^{k+1}}\le
Ch^k\|P\|_{H^1(0,T;H^{k+1})}.
\]

Insert the initial error bounds and
\eqref{eq:residual_estimate_full} into
\eqref{eq:error_stability_with_residuals_full}.
Then combine the resulting estimate with
\[
e_u^n=\eta_u^n+\theta_u^n,
\qquad
e_p^n=\eta_p^n+\theta_p^n,
\]
the projection estimates embedding.
This proves \eqref{eq:error_estimate_full}.
\end{proof}

\begin{remark}
The pressure estimate in the maximum-in-time $V_p$ norm
uses the coercivity of $\mathcal M_p$ on $V_p$.
For this reason, the initial pressure belongs to $V_p$,
and the temporal consistency estimate involves
$\|\partial_{tt}P\|_{L^2(0,T;V_p)}.$

If the storage terms are truncated so that
$G^{pp}=\bar G^{pp}=D^{pp}=0$,
the storage form satisfies
\[
\mathcal M_p(P,P)=(H^{pp}P,P)
\ge \alpha_H\|P\|_{L^2}^2,
\]
but is not generally coercive in $V_p$, and the storage form generally controls only the $L^2$ norm.
Thus, the basic energy estimate controls
$P$ in $L^\infty(0,T;L^2)$, while the diffusion form
still controls $P$ in $L^2(0,T;V_p)$.
The corresponding energy-based error estimate uses
$\max_n\|e_p^n\|_{L^2}^2$ in place of
$\max_n\|e_p^n\|_{V_p}^2$, while retaining
$\tau\sum_n\|e_p^{n+1}\|_{V_p}^2$.
\end{remark}

\section{Numerical results}
\label{sec:numerics}

Having established the well-posedness, stability, and convergence in Sections~\ref{sec:wellposedness} and~\ref{sec:error}, we now present illustrative numerical results. The main goal of this section is to go beyond the theoretical results and to demonstrate the accuracy of the multicontinuum model for complex multiscale media, where the exact solution is taken from the microscale. For this purpose, we solve the fine-scale Biot model \eqref{eq:biot_model} on a fine grid that resolves all heterogeneities and regard its averaged solution as exact. Specifically, the fine-scale pressure and displacements are averaged over each continuum in every coarse block and compared with the corresponding macroscopic variables $P_i$ and $U_i$ in the relative $L^2$ norm. Thus, the obtained errors reflect not only the discretization error but also the modeling error of the multicontinuum approach.

\begin{figure}[hbt!]
\centering
\includegraphics[width=0.3\textwidth]{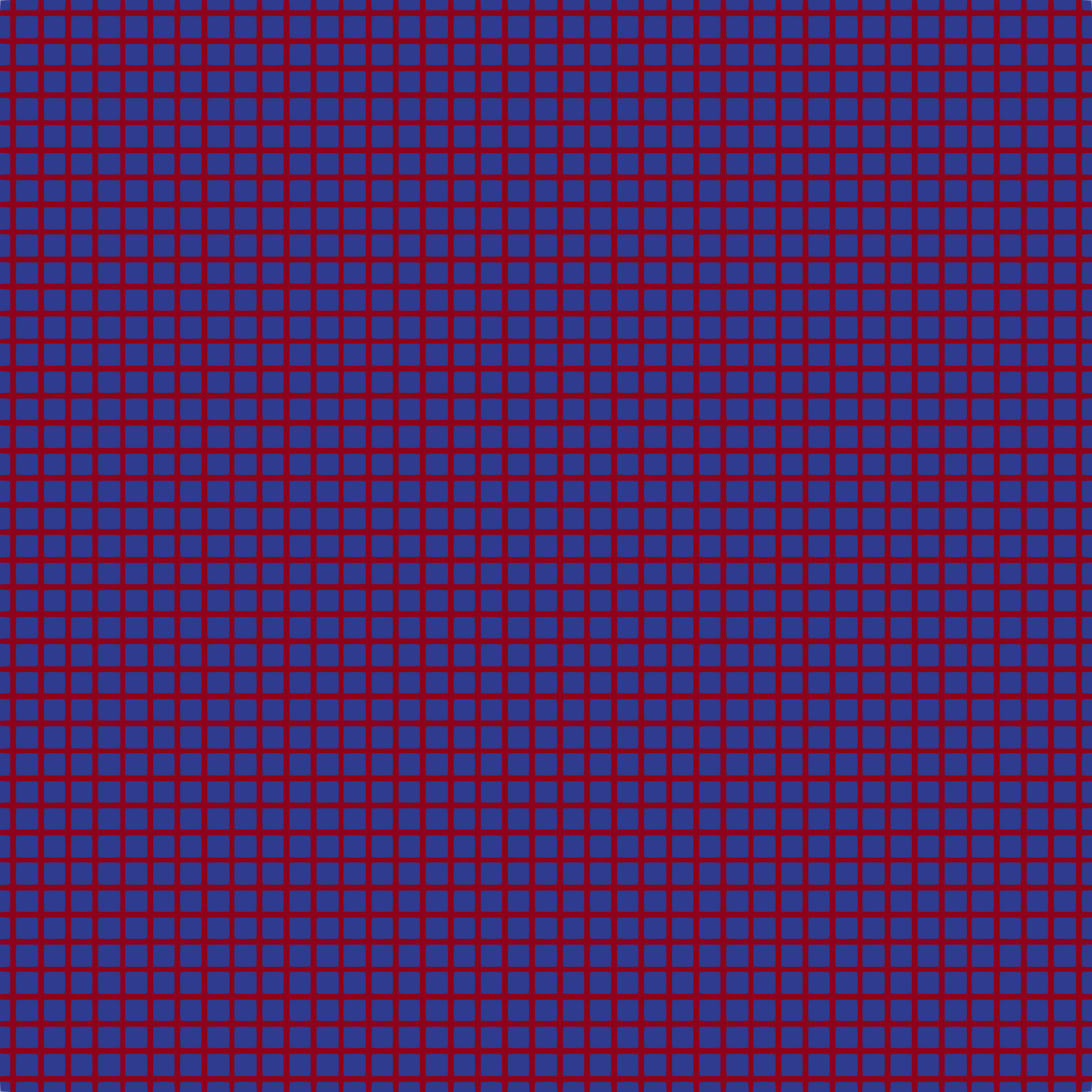}
\includegraphics[width=0.3\textwidth]{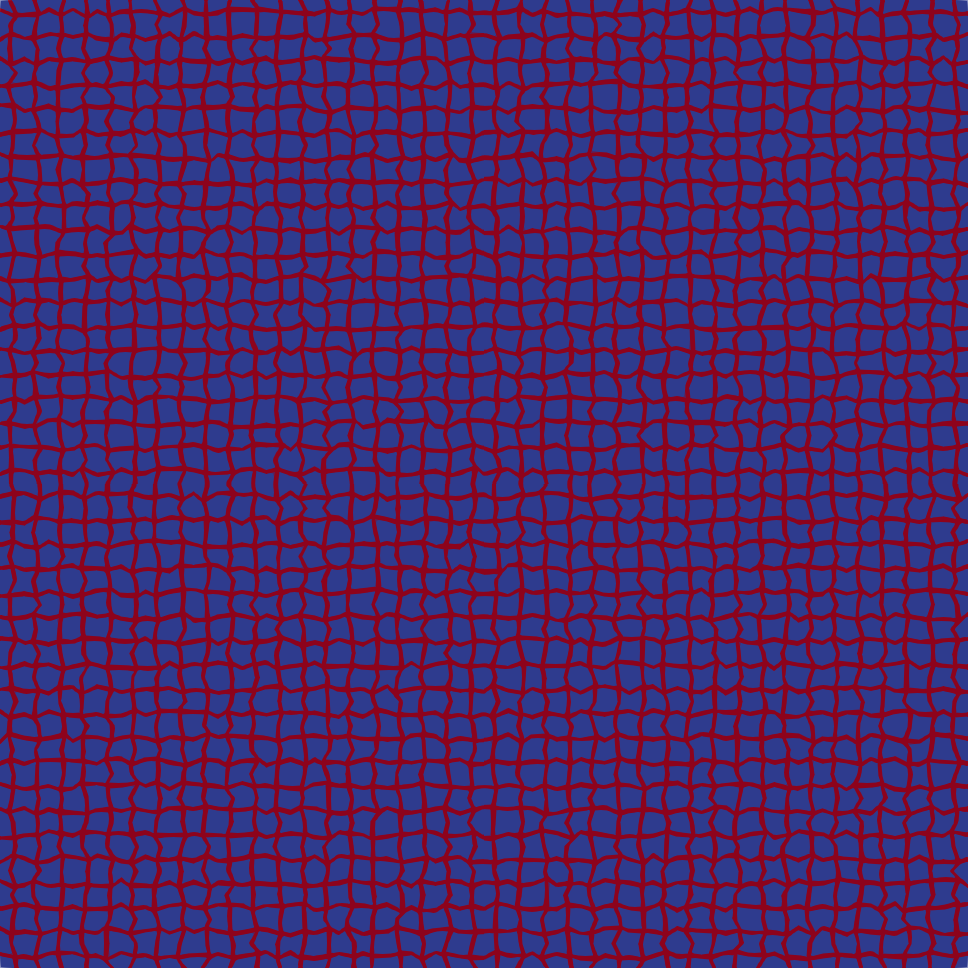}
\caption{Microstructures ($\Omega_1$ is blue, $\Omega_2$ is red). Left: periodic. Right: non-periodic.}
\label{fig:microstructure}
\end{figure}

We consider the computational domain $\Omega=(0,1)\times(0,1)$ occupied by a high-contrast heterogeneous medium. The domain is divided into two non-overlapping subregions, $\Omega=\Omega_1\cup\Omega_2$, which define the two continua ($N_c=2$). In our experiments, we consider both periodic and non-periodic microstructures presented in  Figure~\ref{fig:microstructure}, where $\Omega_1$ is blue and $\Omega_2$ is red.

For our numerical test, we set the following heterogeneous Lam\'e parameters $\lambda$ $[\text{Pa}]$ and $\mu$ $[\text{Pa}]$ and the permeability coefficient $\kappa$ [$\text{m}^2$ / (\text{Pa} \text{s})]
\begin{equation}\label{eq:material_parameters}
\lambda = \begin{cases}
    10^{5}, \quad x \in \Omega_1,\\
    10^{9}, \quad x \in \Omega_2
\end{cases}, \quad
\mu = \begin{cases}
    10^{5}, \quad x \in \Omega_1,\\
    10^{9}, \quad x \in \Omega_2
\end{cases}, \quad
\kappa = \begin{cases}
    10^{-10}, \quad x \in \Omega_1,\\
    10^{-6}, \quad x \in \Omega_2
\end{cases}.
\end{equation}

We impose the boundary conditions
\begin{equation}\label{eq:bcs_ex}
\begin{split}
&u_1 = 0, \quad (\sigma \cdot \nu)_2 = 0, \quad x \in \Gamma_R,\\
&u_2 = 0, \quad (\sigma \cdot \nu)_1 = 0, \quad x \in \Gamma_B,\\
&\sigma \cdot \nu = 0, \quad x \in \Gamma_L \cup \Gamma_T,\\
&-\kappa \frac{\partial p}{\partial \nu} = 0, \quad x \in \partial \Omega,
\end{split}
\end{equation}
where $\nu$ is a unit outward normal vector, and $\Gamma_L$, $\Gamma_R$, $\Gamma_B$, and $\Gamma_T$ denote the left, right, bottom, and top boundaries, respectively, such that $\partial \Omega = \Gamma_L \cup \Gamma_R \cup \Gamma_B \cup \Gamma_T$. Although the analysis in Sections~\ref{sec:wellposedness} and~\ref{sec:error} was presented for homogeneous Dirichlet conditions for simplicity, similar arguments apply to these boundary conditions with minor modifications.

The body force $ f=(f_1,f_2)$ $[\text{N}/\text{m}^3]$ and the fluid source $g$ $[\text{s}^{-1}]$ are given by
\begin{equation}\label{eq:sources}
\begin{gathered}
f_{1} = -10^7\sin(2\pi x_1)\sin(\pi x_2), \quad
f_{2} = 10^7\sin(\pi x_1)\sin(\pi x_2),\\
g = \exp\left(-40\left[(x_1-0.5)^2+(x_2-0.5)^2\right]\right).
\end{gathered}
\end{equation}
The final time is $T=5$ $\text{s}$, and the initial pressure is $p_0=10^6$ $\text{Pa}$. The initial displacement is $u_0 = (0, 0)$ $\text{m}$.

For the periodic microstructure, we employ uniform fine grid consisting of 160801 vertices and 320000 triangular
cells. We construct an unstructured fine grid consisting of 85085 vertices and 169172 triangular cells for the non-periodic microstructure. The fine-scale Biot model \eqref{eq:biot_model} is solved using linear finite elements. For the multicontinuum model, we use two uniform coarse grids: $10\times10$ (121 vertices and 100 square cells) and $20\times20$ (441 vertices and 400 square cells). Each coarse block is taken as the RVE, and the oversampled region is obtained by extending the RVE by $l$ layers of coarse blocks, with $l=5$ for the $10\times10$ grid and $l=6$ for the $20\times20$ grid. The macroscopic equations are discretized using bilinear finite elements on the coarse grid. In time, both the fine-scale and multicontinuum problems are solved with the backward Euler method using 50 uniform time steps. The numerical implementation is based on the FEniCS package \cite{logg2012automated}, and the fine grid is generated using Gmsh \cite{geuzaine2009gmsh}.

The coarse-grid results are based on the multicontinuum model \eqref{eq:multicontinum}. It should be noted that we solve split constraint cell problems in the oversampled RVE for the periodic microstructure. For the non-periodic microstructure, we employ coupled approach by solving coupled cell problems (following \cite{Ammosov2025}).

To compare the multicontinuum and reference solutions, we use the following
relative $L^2$ errors for pressure and displacements at the final time
\[
\begin{gathered}
e_p^{(i)}=\sqrt{\frac{\sum_K\left|\frac{1}{|K|}\int_K P_i\,dx-\frac{1}{|K\cap\Omega_i|}\int_{K\cap\Omega_i}p\,dx\right|^2}{\sum_K\left|\frac{1}{|K\cap\Omega_i|}\int_{K\cap\Omega_i}p\,dx\right|^2}},\\[0.5em]
e_u^{(i)}=\sqrt{\frac{\sum_K\left|\frac{1}{|K|}\int_K U_i\,dx-\frac{1}{|K\cap\Omega_i|}\int_{K\cap\Omega_i} u\,dx\right|^2}{\sum_K\left|\frac{1}{|K\cap\Omega_i|}\int_{K\cap\Omega_i} u\,dx\right|^2}},
\end{gathered}
\]
where $i=1,2$ denotes the continuum number, $K$ is the coarse block taken as
the RVE, and $( u,p)$ is the fine-scale reference solution.

\begin{figure}[!htbp]
\centering
\begin{subfigure}{\textwidth}
\centering
\includegraphics[height=0.24\textheight]{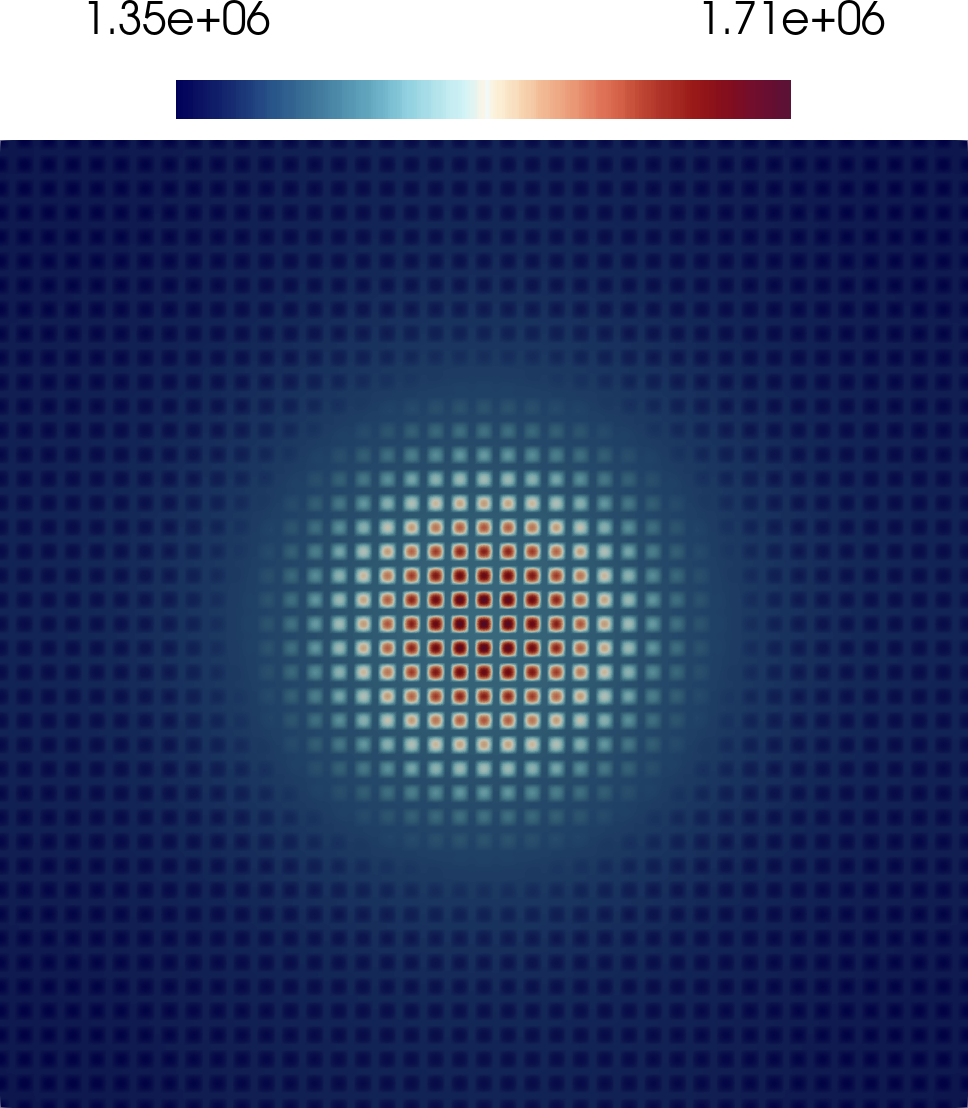}
\includegraphics[height=0.24\textheight]{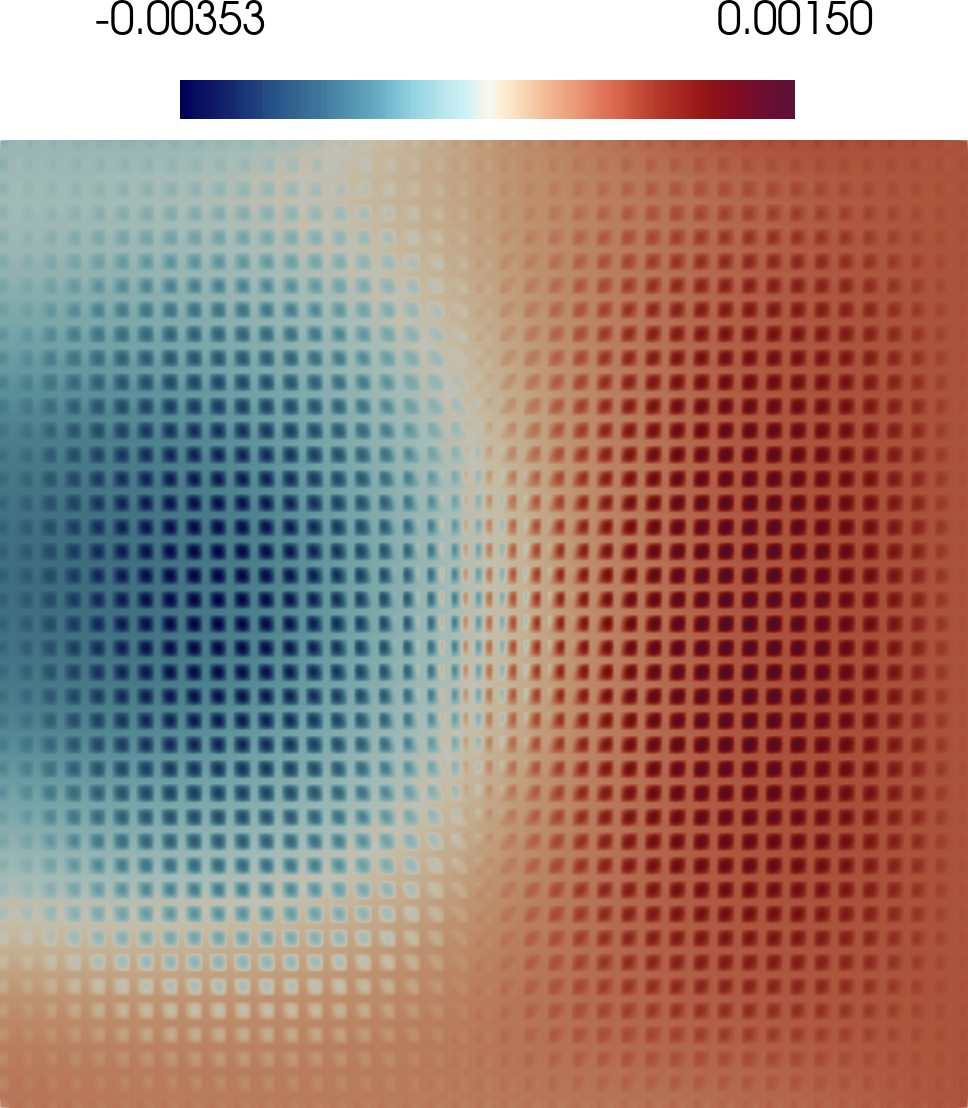}
\includegraphics[height=0.24\textheight]{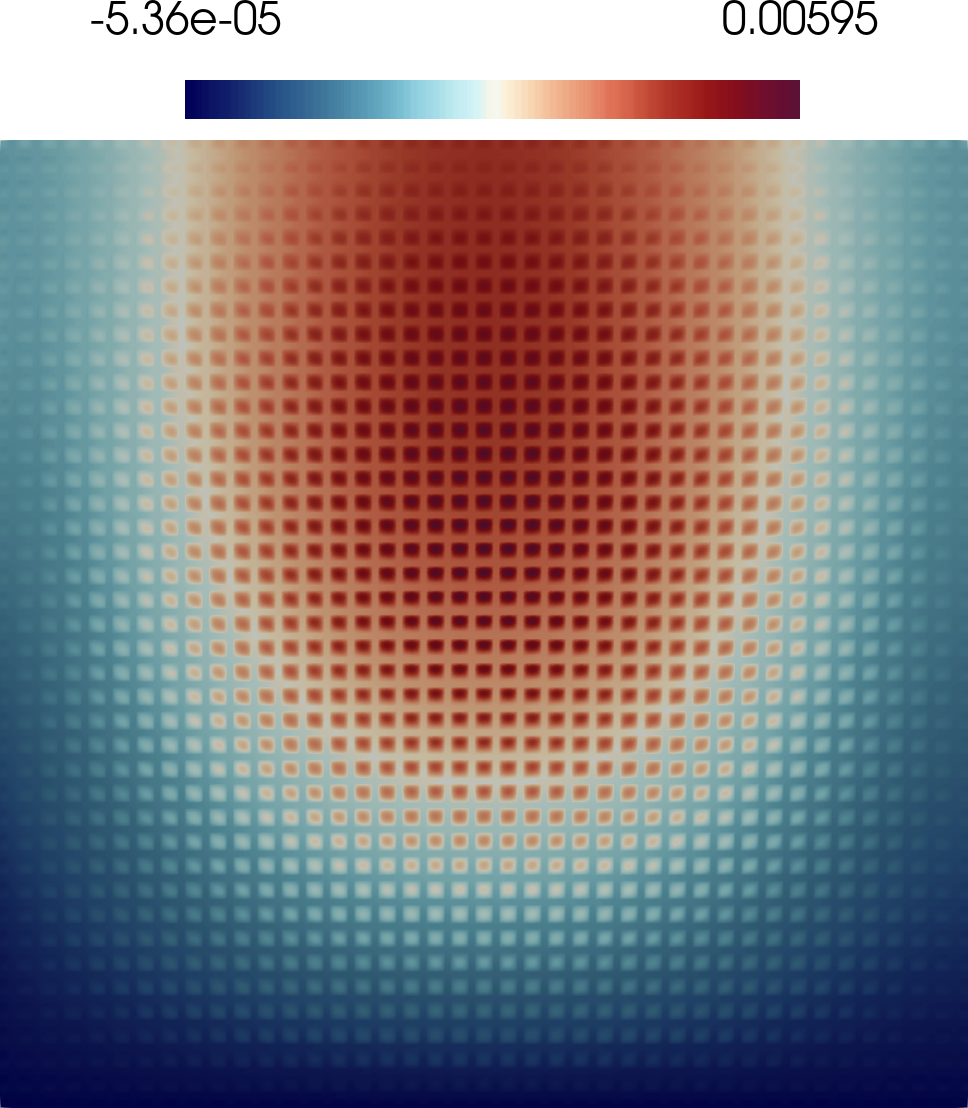}
\caption{Fine-grid solution}
\label{fig:results_periodic_fine}
\end{subfigure}
\begin{subfigure}{\textwidth}
\centering
\includegraphics[height=0.24\textheight]{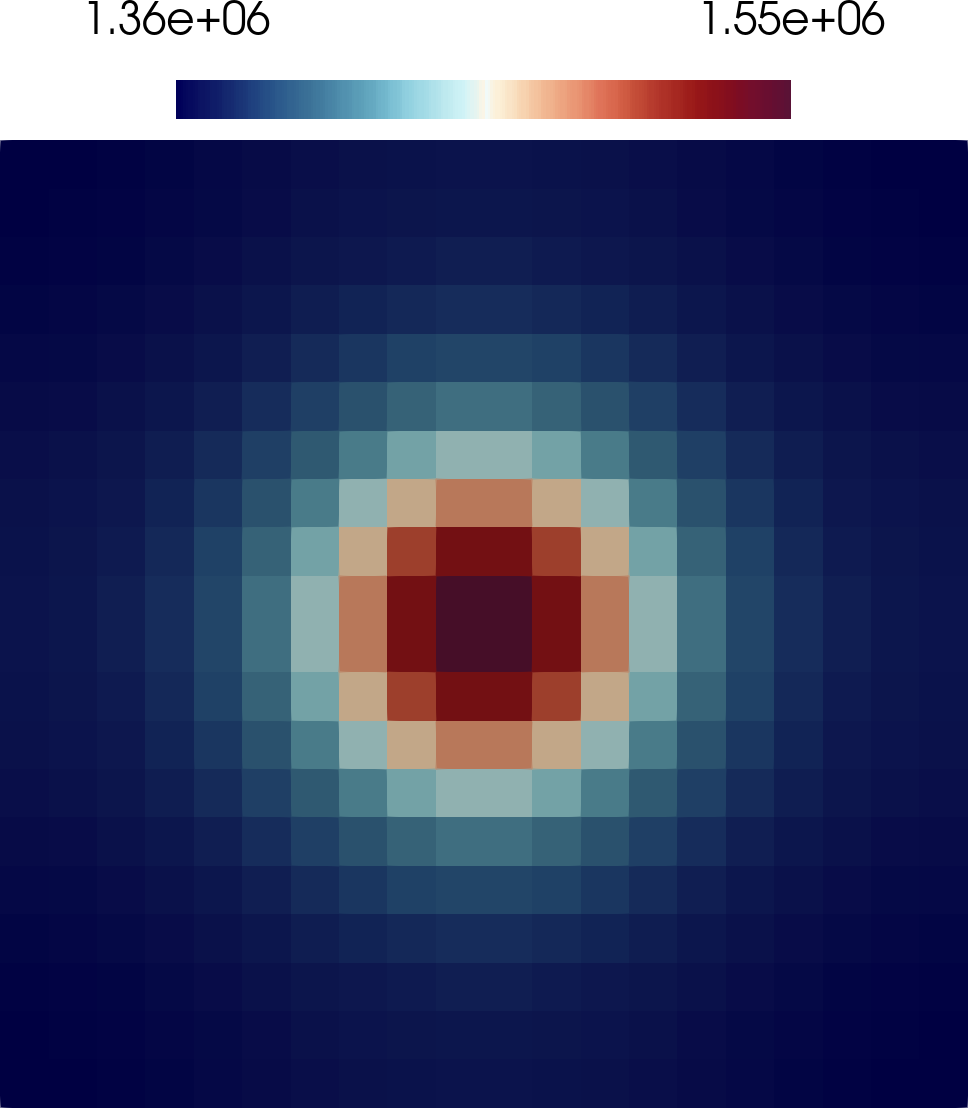}
\includegraphics[height=0.24\textheight]{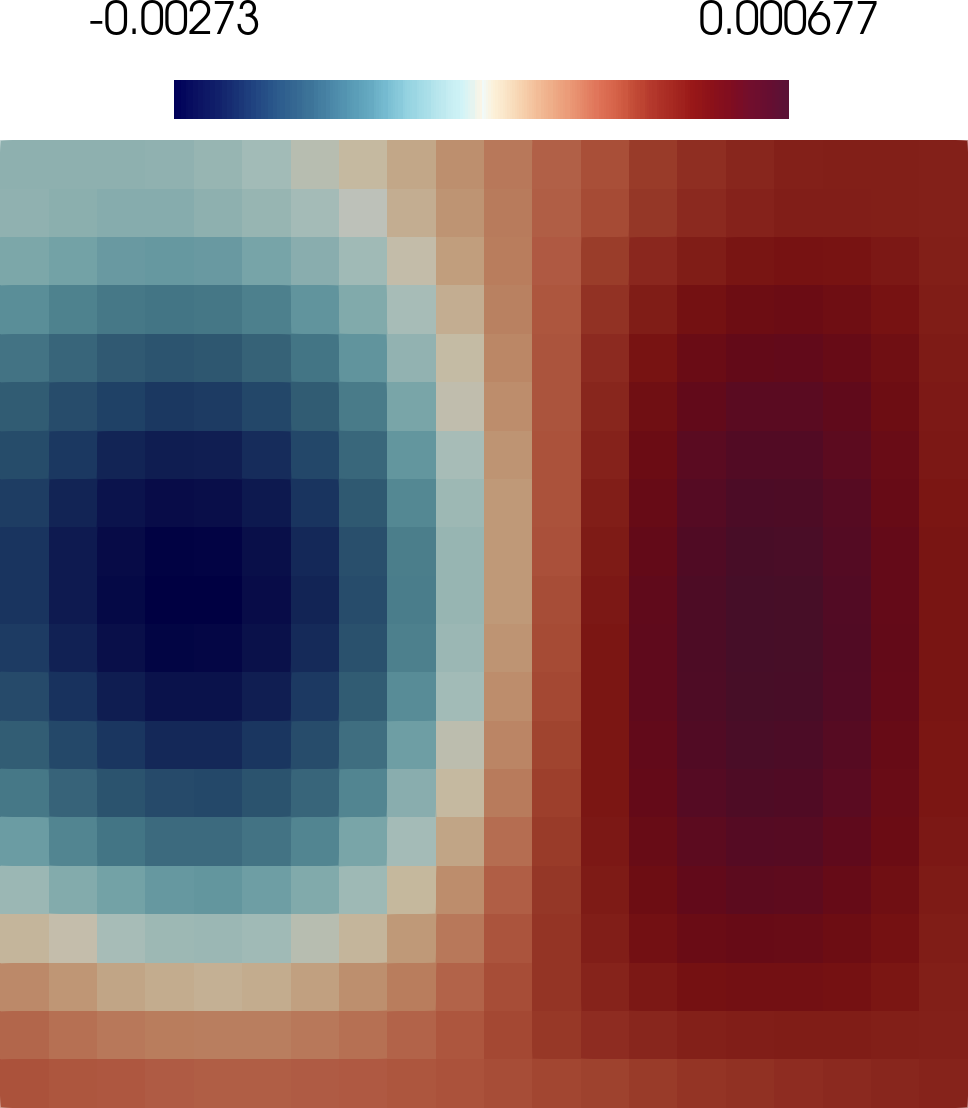}
\includegraphics[height=0.24\textheight]{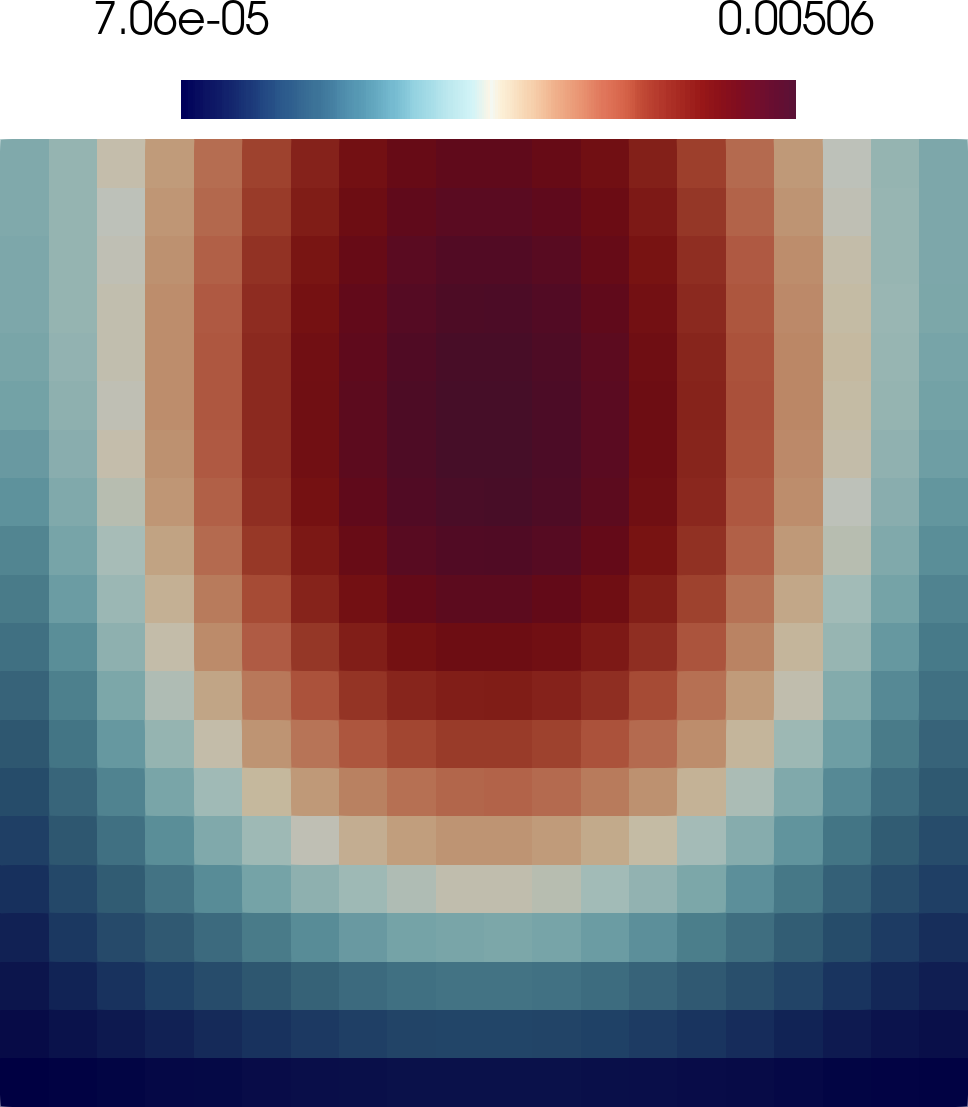}
\caption{Average solution in $\Omega_1$ on the coarse grid $20 \times 20$}
\label{fig:results_periodic_coarse_continuum_1}
\end{subfigure}
\begin{subfigure}{\textwidth}
\centering
\includegraphics[height=0.24\textheight]{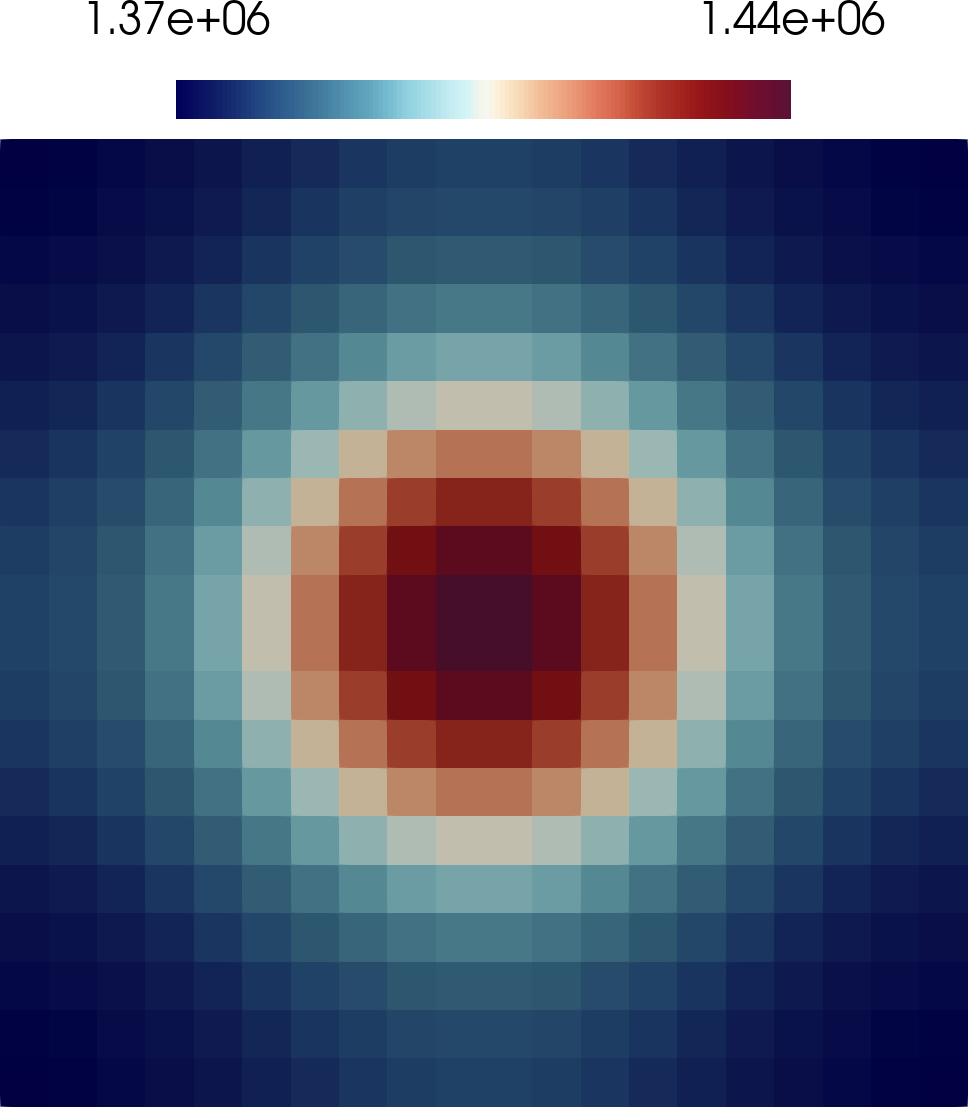}
\includegraphics[height=0.24\textheight]{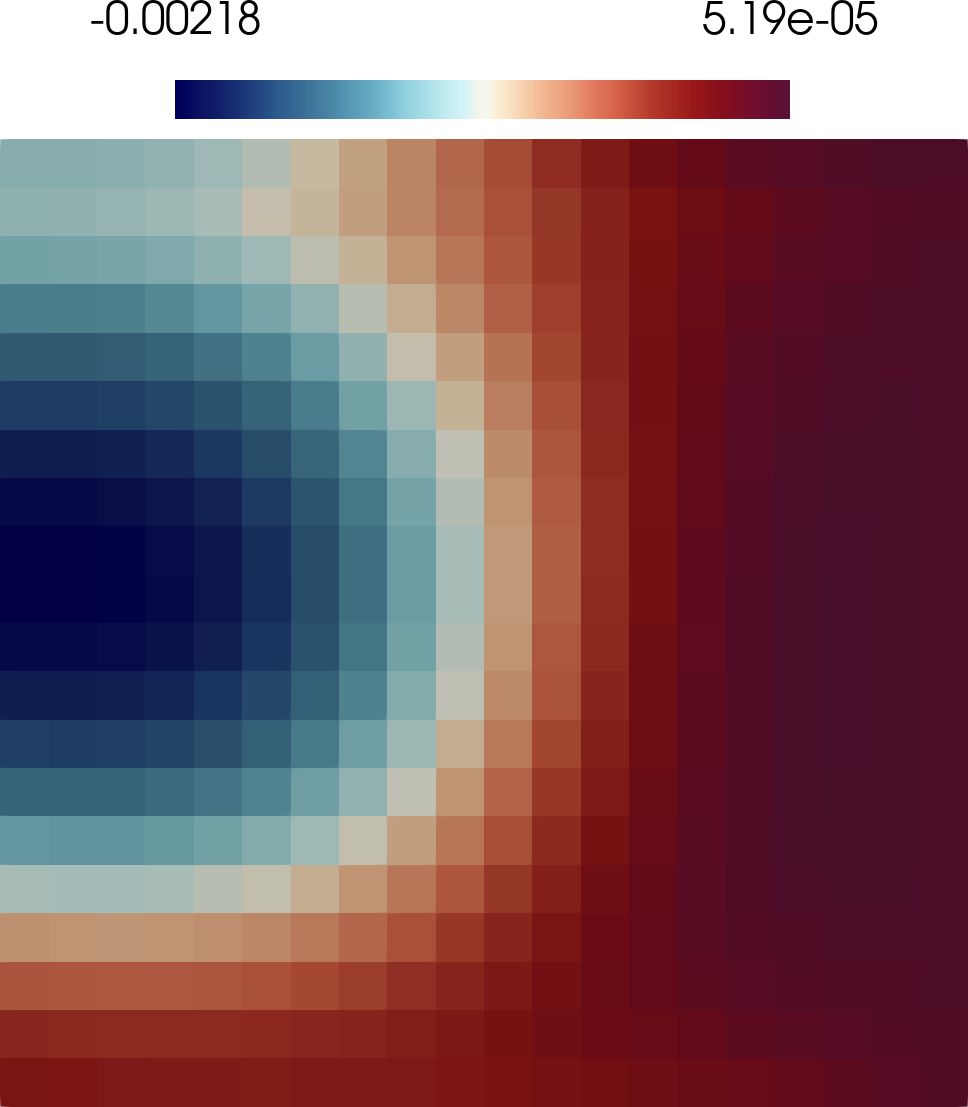}
\includegraphics[height=0.24\textheight]{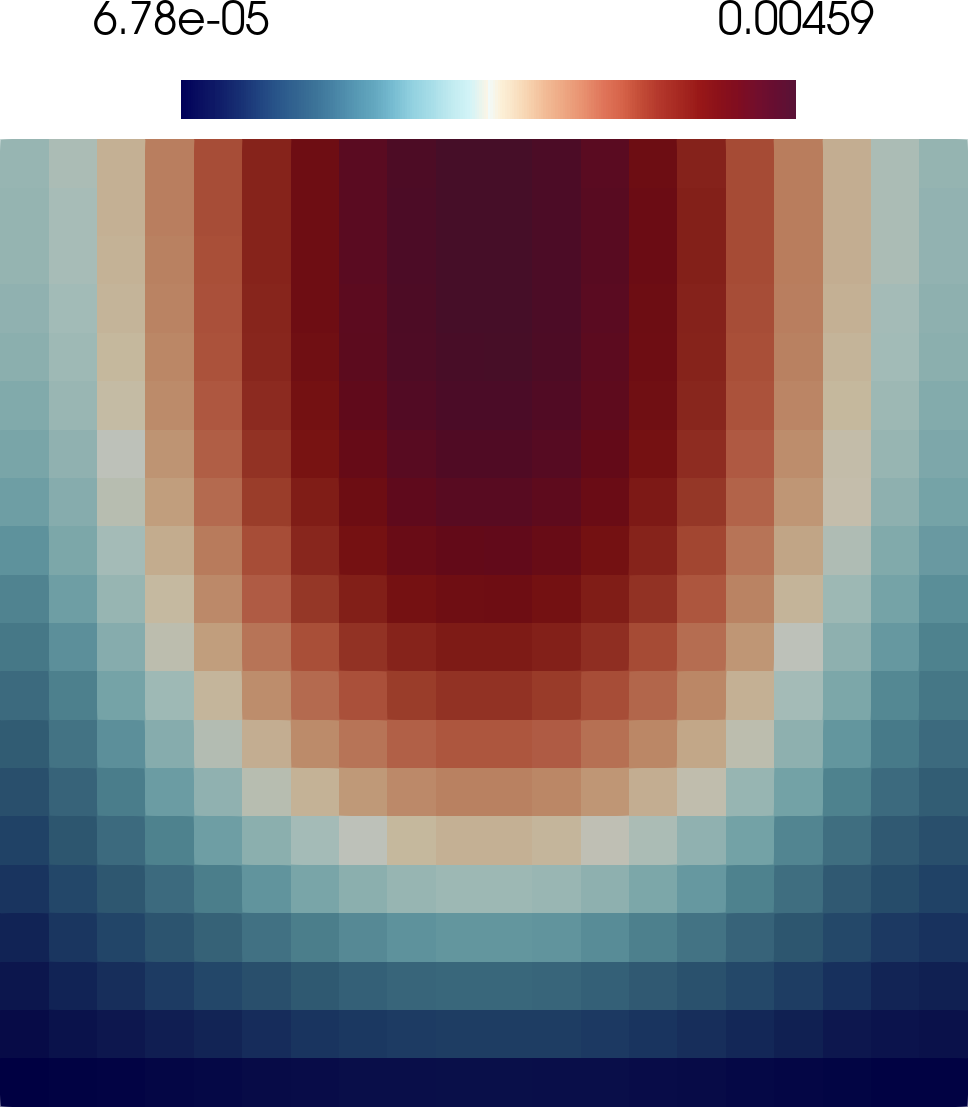}
\caption{Average solution in $\Omega_2$ on the coarse grid $20 \times 20$}
\label{fig:results_periodic_coarse_continuum_2}
\end{subfigure}
\caption{Distributions of reference pressure $[\text{Pa}]$ and displacements $[\text{m}]$ in $x_1$ and $x_2$ directions (from left to right) at the final time. Periodic microstructure.}
\label{fig:results_periodic}
\end{figure}

\begin{figure}[!htbp]
\centering
\begin{subfigure}{\textwidth}
\centering
\includegraphics[height=0.24\textheight]{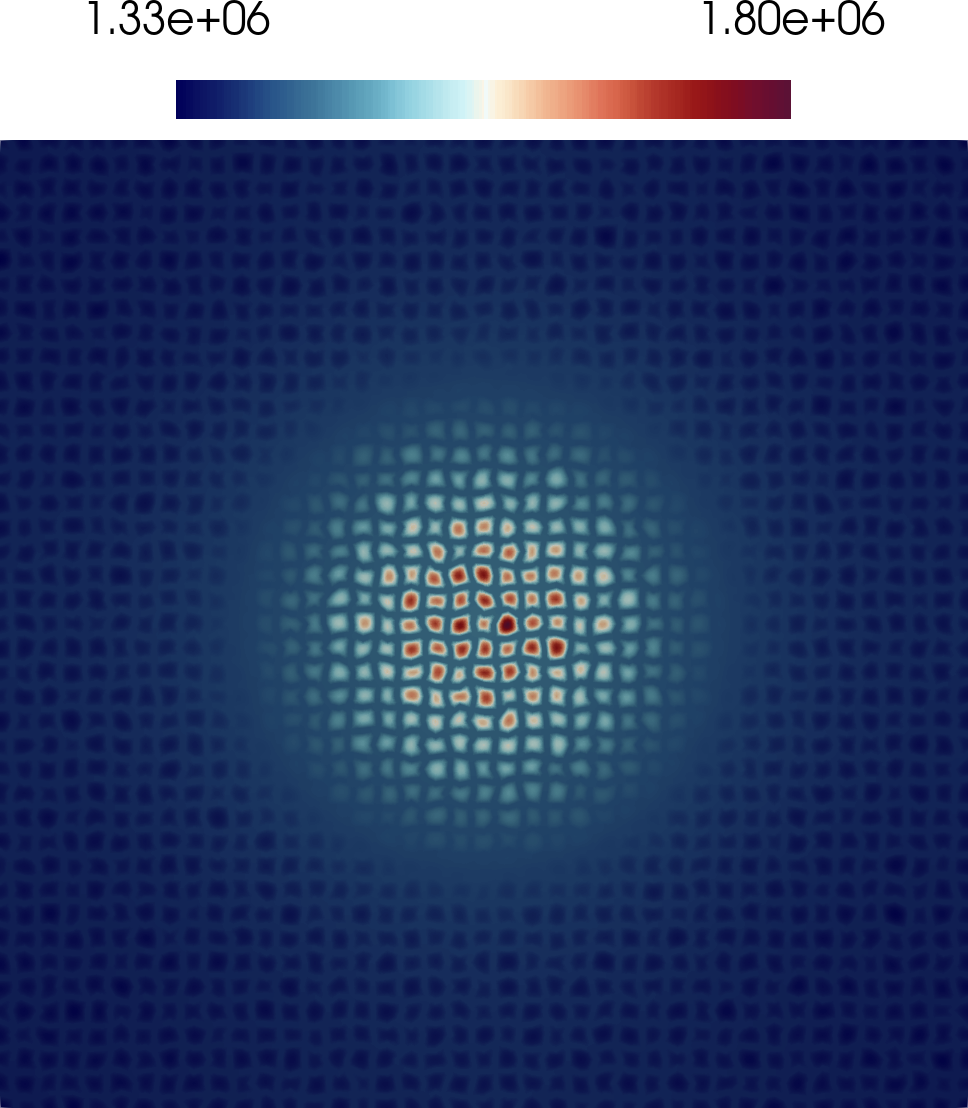}
\includegraphics[height=0.24\textheight]{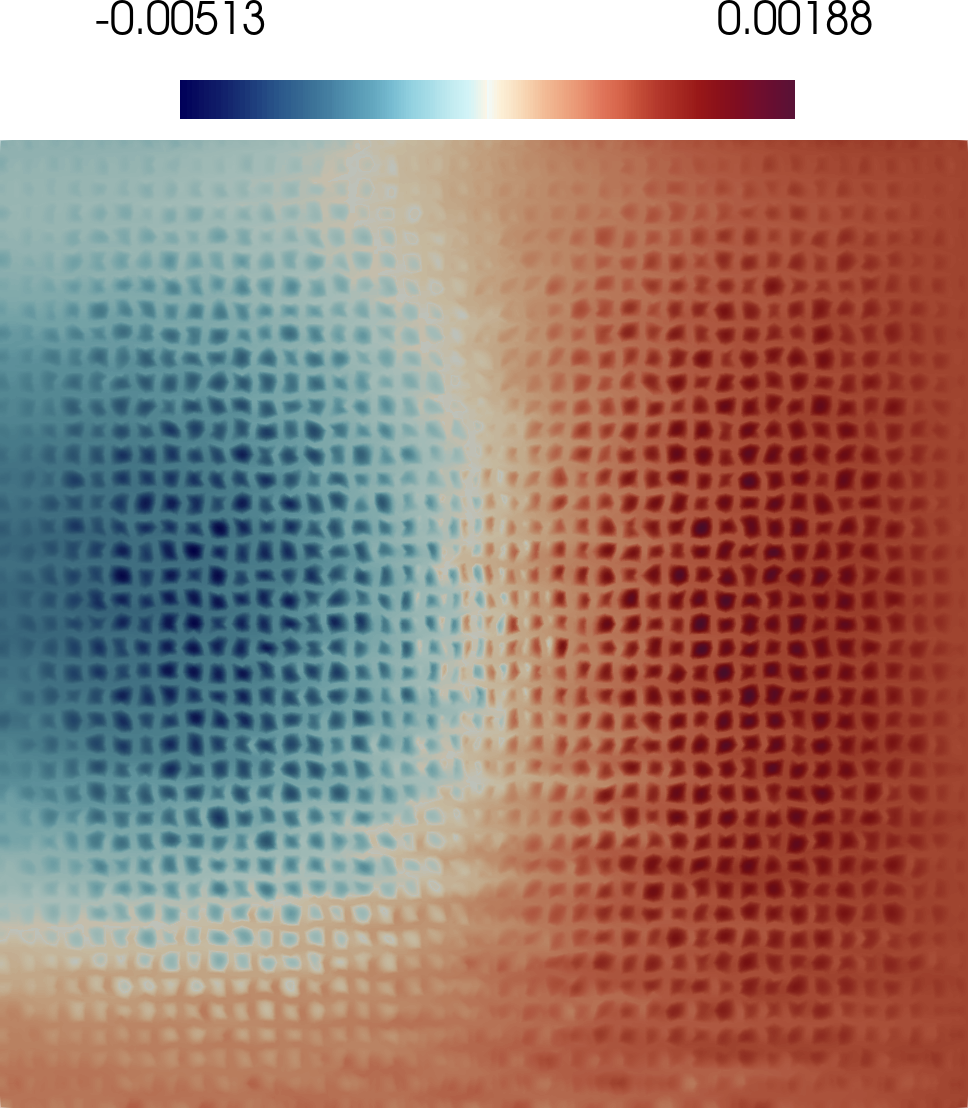}
\includegraphics[height=0.24\textheight]{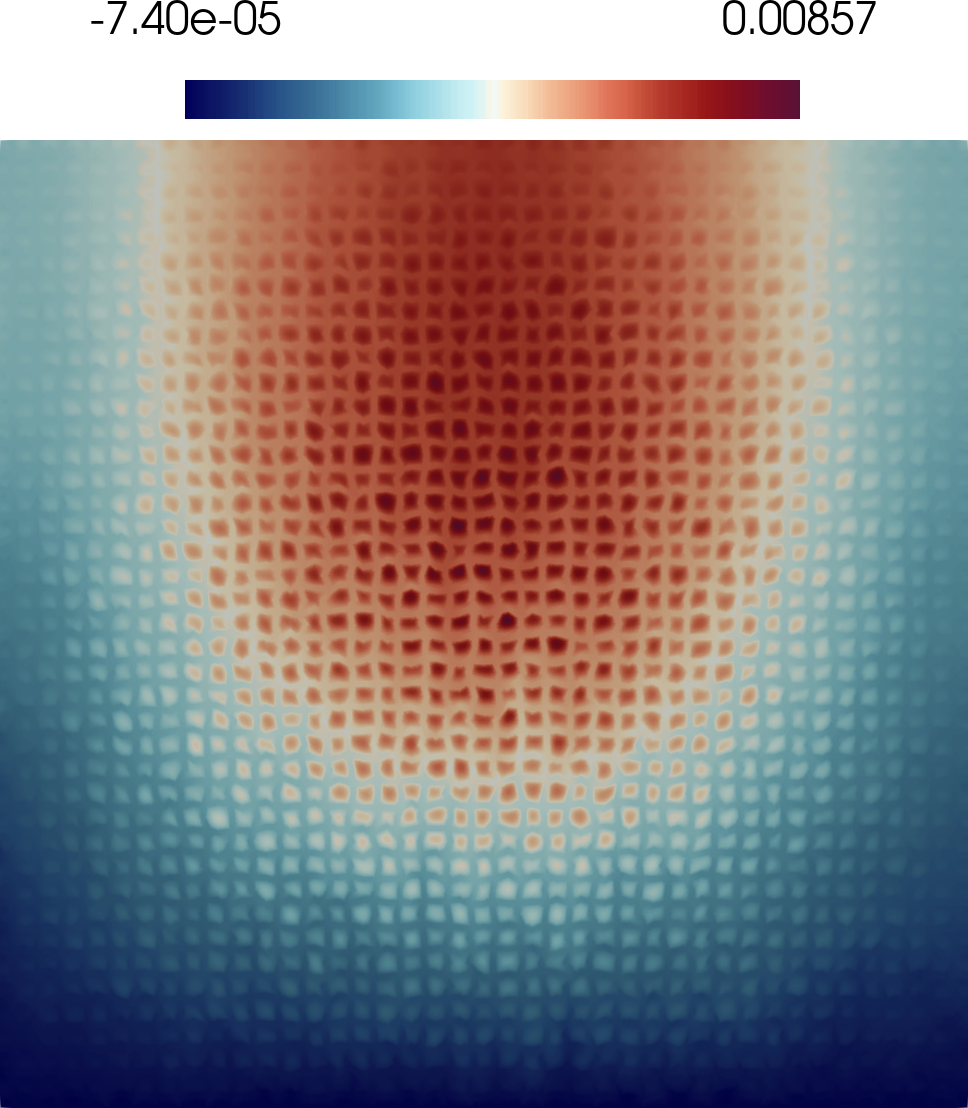}
\caption{Fine-grid solution}
\label{fig:results_non_periodic_fine}
\end{subfigure}
\begin{subfigure}{\textwidth}
\centering
\includegraphics[height=0.24\textheight]{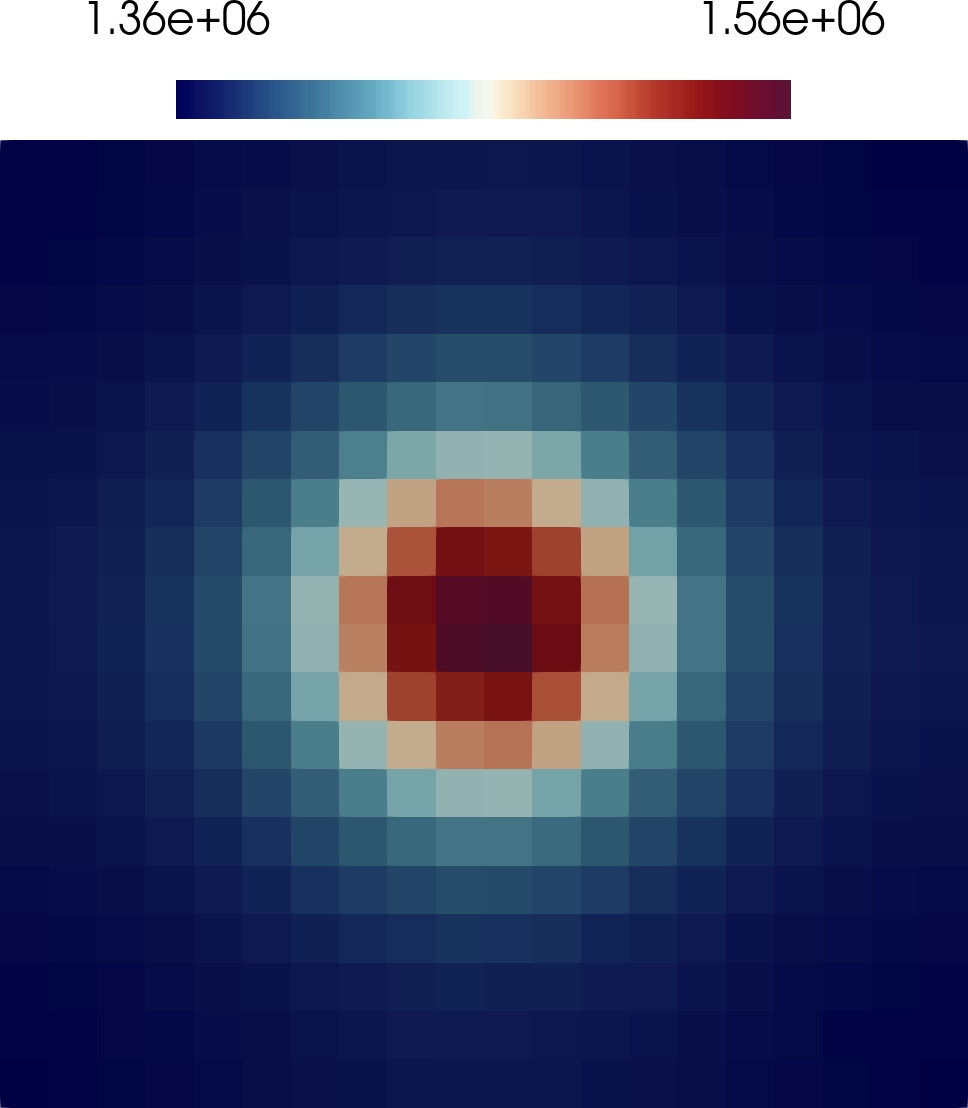}
\includegraphics[height=0.24\textheight]{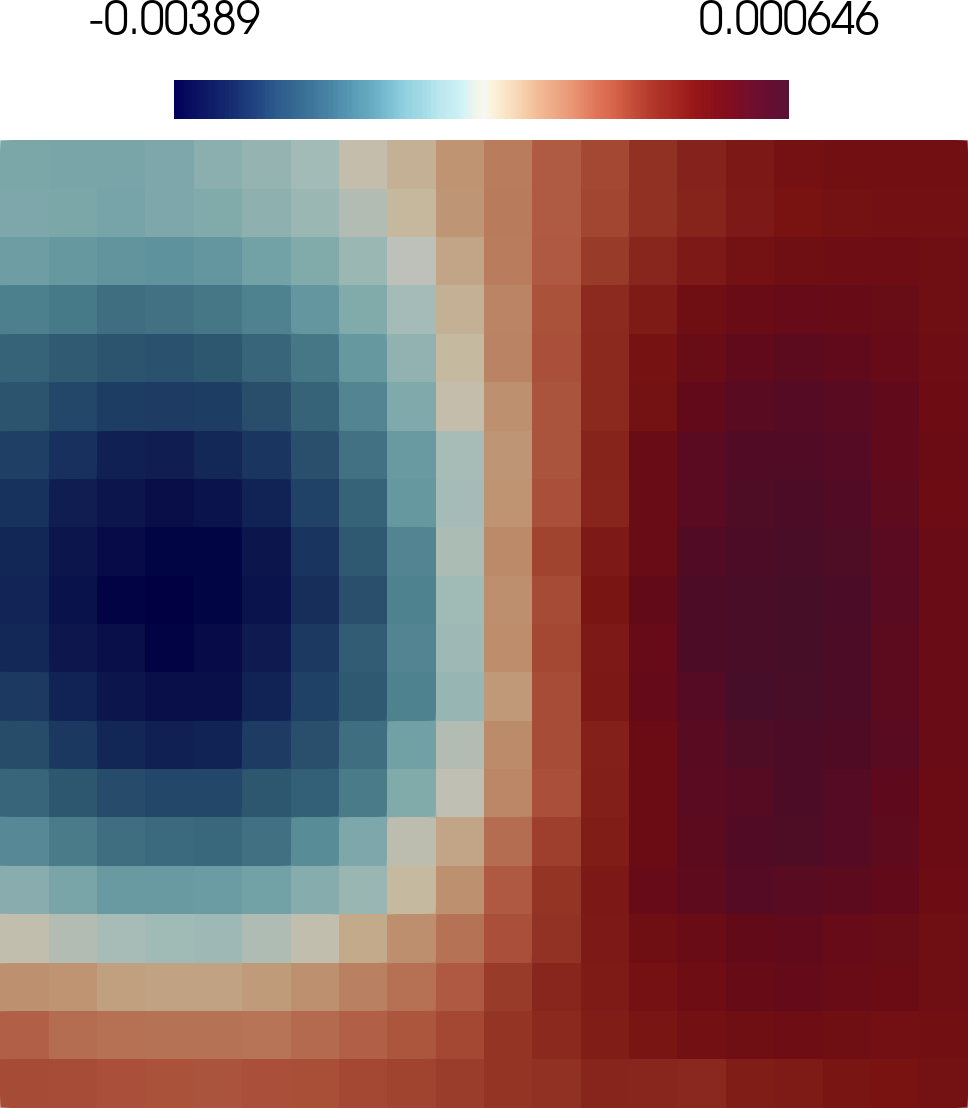}
\includegraphics[height=0.24\textheight]{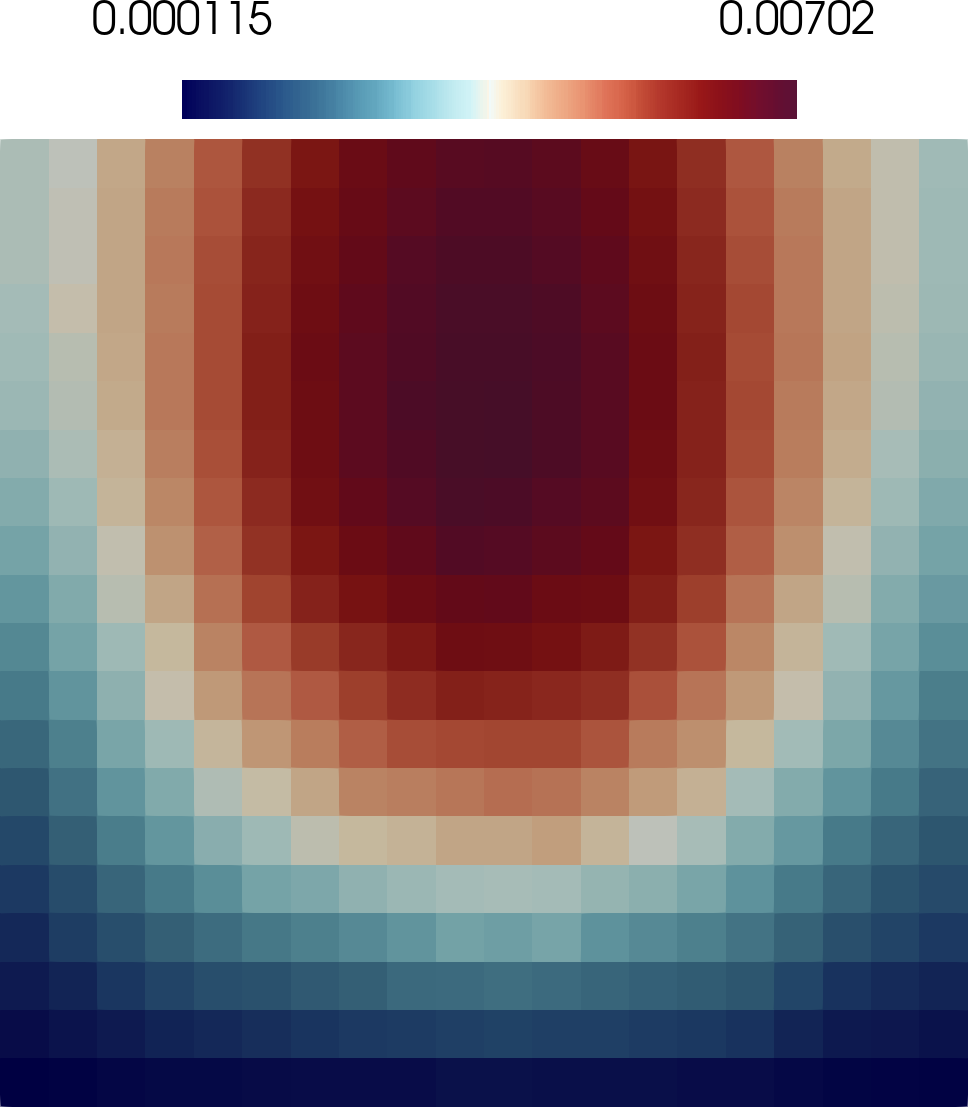}
\caption{Average solution in $\Omega_1$ on the coarse grid $20 \times 20$}
\label{fig:results_non_periodic_coarse_continuum_1}
\end{subfigure}
\begin{subfigure}{\textwidth}
\centering
\includegraphics[height=0.24\textheight]{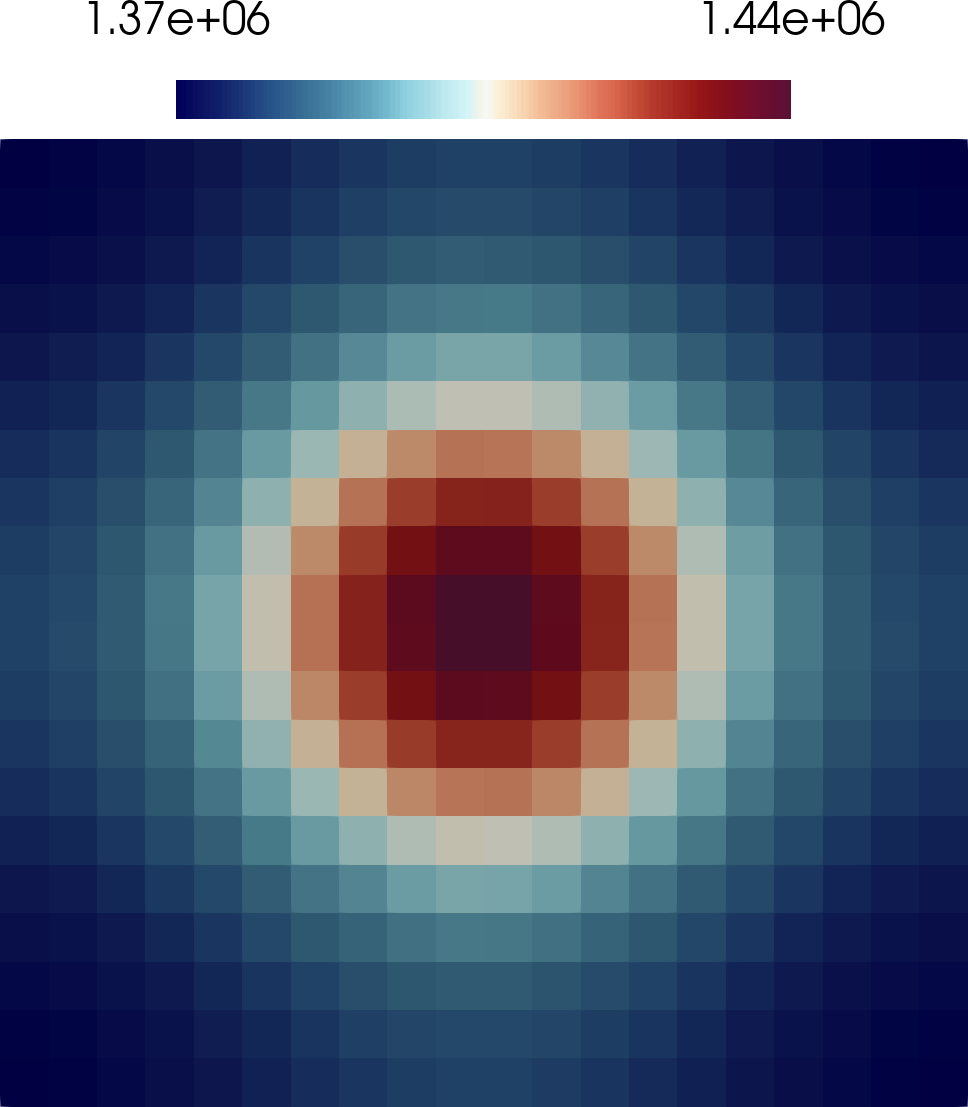}
\includegraphics[height=0.24\textheight]{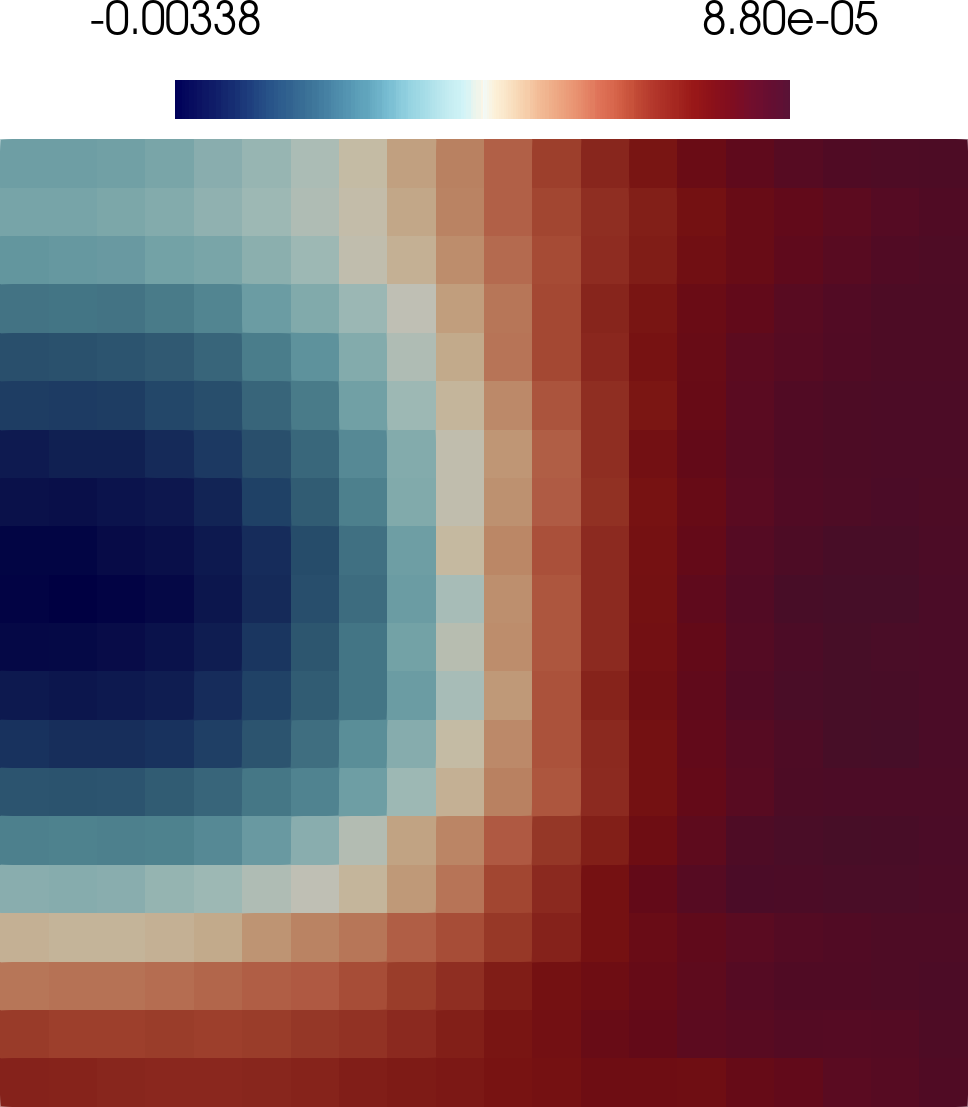}
\includegraphics[height=0.24\textheight]{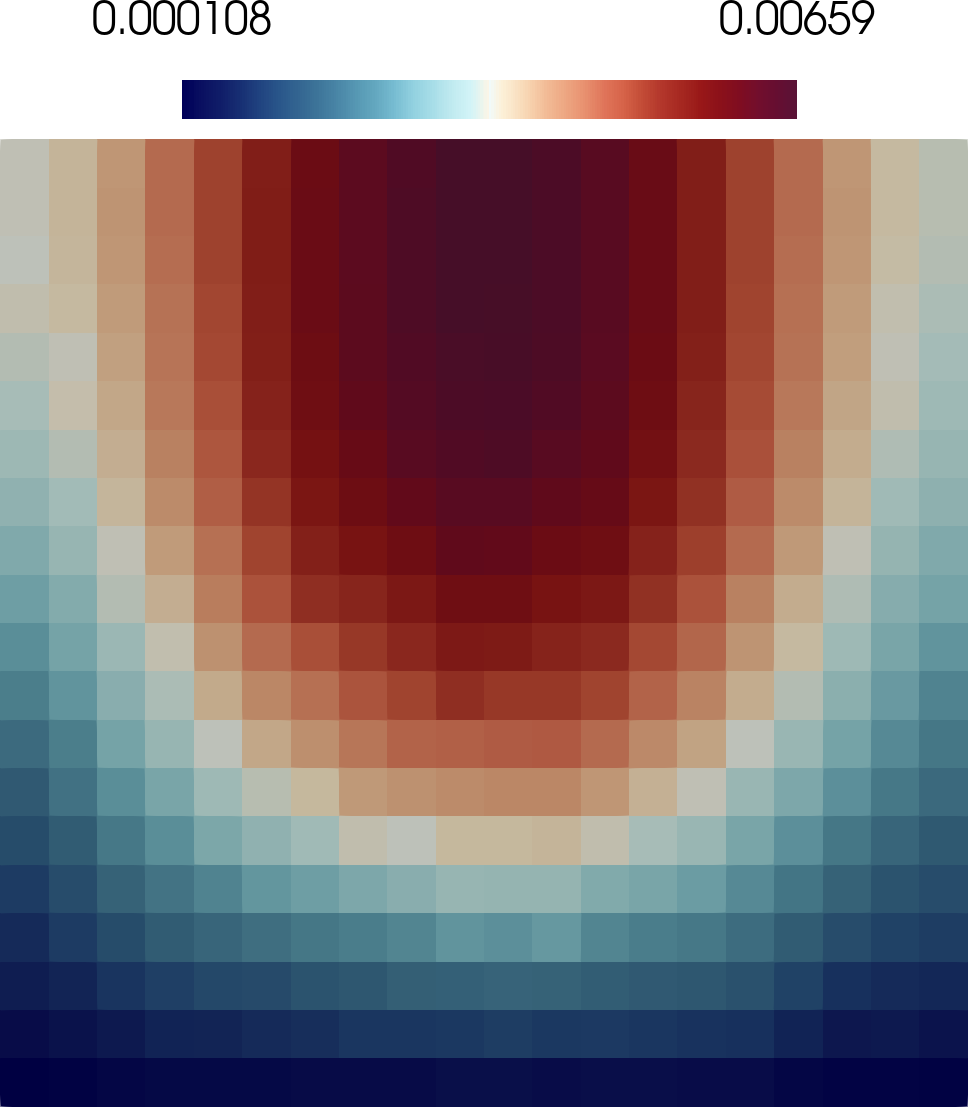}
\caption{Average solution in $\Omega_2$ on the coarse grid $20 \times 20$}
\label{fig:results_non_periodic_coarse_continuum_2}
\end{subfigure}
\caption{Distributions of reference pressure $[\text{Pa}]$ and displacements $[\text{m}]$ in $x_1$ and $x_2$ directions (from left to right) at the final time. Non-periodic microstructure.}
\label{fig:results_non_periodic}
\end{figure}

Figures~\ref{fig:results_periodic} and~\ref{fig:results_non_periodic} show the reference solutions at the final time for the periodic and non-periodic microstructures, respectively. In both cases, the fine-grid fields exhibit small-scale oscillations induced by the heterogeneous microstructure. The pressure builds up around the source located at the center of the domain and remains elevated in the low-permeability subregion $\Omega_1$. The horizontal displacement changes sign across the domain, indicating that the medium is pulled apart in the $x_1$ direction, whereas the vertical displacement is positive and peaks in the upper central part, which reflects the upward deformation. On the coarse grid $20\times20$, the continuum averages retain these features. The averaged pressure is more widely distributed in $\Omega_2$, where the permeability is higher, and is more concentrated and larger near the source in $\Omega_1$. Both averaged displacement components have larger magnitudes in the softer subregion $\Omega_1$. The non-periodic microstructure yields larger pressure and displacement values than the periodic one, but the overall behavior of the solutions is similar.

Figures~\ref{fig:results_periodic_coarse_continuum_1}--\ref{fig:results_periodic_coarse_continuum_2} and~\ref{fig:results_non_periodic_coarse_continuum_1}--\ref{fig:results_non_periodic_coarse_continuum_2} present the reference solutions averaged over each continuum on the coarse grid $20\times20$. The averaged fields preserve the main features of the fine-grid solutions while removing the small-scale oscillations. The averaged pressure in $\Omega_1$ is concentrated near the source and attains larger values, whereas in $\Omega_2$ it is more widely distributed owing to the higher permeability. For both displacement components, the averages in the softer subregion $\Omega_1$ have larger magnitudes than those in $\Omega_2$. Thus, the two continua exhibit distinct responses, which motivates the use of separate macroscopic variables for each of them.

Table~\ref{tab:errors} presents the relative $L^2$ errors at the final time for both microstructures and coarse grids. In all cases, the errors are small: below $0.3\%$ for the pressure and below $5\%$ for the displacements. Moreover, the errors decrease as the coarse grid is refined. For the non-periodic microstructure, the displacement errors are larger and decrease more slowly than in the periodic case. Overall, the results demonstrate that the multicontinuum model provides an accurate approximation of the averaged microscale solution for both periodic and non-periodic media.

\begin{table}[hbt!]
\centering
\begin{tabular}{c|c|c|c|c}
Coarse grid & $e^{(1)}_p$ & $e^{(2)}_p$ & $e^{(1)}_u$ & $e^{(2)}_u$\\ \hline
\multicolumn{5}{c}{Periodic microstructure} \\ \hline
$10 \times 10$ & 2.24e-03 & 6.46e-04 & 1.89e-02 & 1.60e-02 \\
$20 \times 20$ & 6.95e-04 & 1.83e-04 & 5.45e-03 & 3.99e-03 \\ \hline
\multicolumn{5}{c}{Non-periodic microstructure} \\ \hline
$10 \times 10$ & 2.30e-03 & 7.65e-04 & 4.00e-02 & 4.19e-02 \\
$20 \times 20$ & 8.94e-04 & 2.26e-04 & 2.80e-02 & 2.99e-02 \\
\end{tabular}
\caption{Relative $L_2$ errors for different microstructures and coarse grids.}
\label{tab:errors}
\end{table}

In summary, the numerical results show that the multicontinuum model accurately reproduces the averaged microscale solution for both periodic and non-periodic high-contrast media, with errors that decrease as the coarse grid is refined. These results complement the theoretical analysis by demonstrating the accuracy of the model with respect to the underlying microscale problem.

\section{Conclusion}
\label{sec:conclusions}
In this work, we developed and analyzed a macroscopic multicontinuum model for Biot poroelasticity. Starting from the weak formulation of the classical Biot system, we constructed continuum-dependent Galerkin approximation spaces enriched with first-order multiscale corrections. By projecting the microscopic variational formulation onto these spaces, we systematically derived the corresponding macroscopic multicontinuum equations and their effective coefficients. In particular, the resulting model retains the generalized pressure storage and coupling terms arising from the first-order multiscale approximation, providing a consistent variational framework for describing interactions among multiple continua and the deformable solid matrix.

We investigated the mathematical properties of the full macroscopic multicontinuum system. Under suitable assumptions on the effective coefficients, we established the existence and uniqueness of weak solutions and derived an energy stability estimate. By exploiting the properties of the effective coupling operators, we obtained an energy structure in which the coupling contributions cancel exactly. Based on this structure, we developed a fully discrete Galerkin finite element scheme using the backward Euler method for temporal discretization. The well-posedness and stability of the fully discrete scheme were established, and a priori error estimates were derived through suitable projection operators and discrete energy arguments. Numerical experiments for periodic and non-periodic high-contrast media demonstrate the accuracy of the multicontinuum model with respect to the averaged fine-scale solution.

\bibliographystyle{unsrt}
\bibliography{lit}

\end{document}